\documentclass[twoside,letterpaper,draft,11pt]{amsart}
\usepackage{amsmath,amsthm, amscd, amssymb, amsfonts,enumerate}
\usepackage{latexsym,xspace,mathtools,color}
\usepackage[mathscr]{eucal}
\usepackage[all]{xy}
\usepackage{mathrsfs}
\usepackage[T2A,T1]{fontenc}
\usepackage{caption}
\usepackage{graphicx,fancyhdr}
\usepackage[usenames, dvipsnames]{xcolor}
\usepackage{nccmath, relsize}
\usepackage[normalem]{ulem}

\usepackage[utf8]{inputenc}

\newcommand{\m}{{}^{-1}}
\newcommand{\ta}{\tau}
\newcommand{\bt}{\beta}
\newcommand{\af}{\alpha}

\newcommand{\comment}[1]{}

\newtheorem{teo1}{Theorem}[section]
\newtheorem{def1}[teo1]{Definition}
\newtheorem{lem1}[teo1]{Lemma}
\newtheorem{cor1}[teo1]{Corollary}
\newtheorem{prop1}[teo1]{Proposition}

\newtheorem{exe}[teo1]{Example}
\newtheorem{remark}[teo1]{Remark}

\def\pf{\begin{proof}}
\def\epf{\end{proof}}

\DeclareMathOperator{\G}{G}
\DeclareMathOperator{\cH}{H}
\DeclareMathOperator{\moL}{L}

\newcommand{\Subg}[1]{\mathfrak{L}(#1)}
\newcommand{\wSubg}[1]{\mathfrak{L}_{w}(#1)}

\newcommand{\id}{\operatorname{id}}

\begin{document}

\title[Galois correspondence for groupoid partial actions]{Galois correspondence for group-type partial actions of groupoids}

\author[Bagio, Sant'Ana, Tamusiunas]
{Dirceu Bagio, Alveri Sant'Ana, Tha\'{\i}sa Tamusiunas}

\address{Departamento de Matem\'atica, Universidade Federal de Santa Maria,
	97105-900, Santa Maria, RS, Brazil} \email{bagio@smail.ufsm.br}

\address{Instituto de Matem\'atica e Estat\'istica, Universidade federal de Porto Alegre, 91509-900\\
	Porto Alegre-RS, Brazil}
\email{alveri@ufrgs.br, thaisa.tamusiunas@gmail.com}


\begin{abstract}
Let $\G$ be a finite groupoid and $\alpha=(S_g,\af_g)_{g\in \G}$ a unital partial action of group-type of $\G$ on a commutative ring $S=\oplus_{y\in\G_0}S_y$. We shall prove a Galois correspondence between a class of wide subgroupoids of $\G$ and a class of subrings of $S$. We recover known results for global groupoid actions and we give several examples to illustrate the correspondence.
\end{abstract}

\maketitle

\section{Introduction}
The usual notion of Galois extension over fields was extended for commutative rings by M. Auslander and O. Goldman in \cite{AG}. Some years later, the Galois theory over commutative rings was developed by S. U. Chase, D. K. Harrison and A. Rosenberg in \cite{CHR}. They presented several equivalent conditions for the definition of Galois extension. Among the main results, they proved a Galois correspondence in the context of commutative rings. Precisely, if $R\subset S$ is a Galois extension of commutative rings with Galois group $G$, then there exists a bijective association between the set of subgroups of $G$ and the set of $R$-subalgebras of $S$ which are $G$-strong and $R$-separable.    \vspace{0.07cm}

In the 1990's, R. Exel introduced the notion of partial actions of a group  in the theory of operator algebras, see for instance \cite{EX1} and \cite{EX2}. The same notion in an algebraic context was considered in \cite{DE}. Particularly, it was defined partial actions of groups on rings which is the key to develop a partial Galois theory. So, the Galois theory for partial actions of groups on rings was presented two years later in \cite{DFP} generalizing the results of \cite{CHR}.  \vspace{0.07cm}

On the other hand, in the context of category theory, a groupoid is a small category in which every morphism has inverse. However, a groupoid can be seen as a natural generalization of a group. In fact, a groupoid is a set $\G$ equipped with a set of identities $\G_0\subset \G$ and a binary operation defined partially which is associative and, for each $g\in \G$, there exist $g\m\in \G$ such that $g\m g=s(g)\in \G_0$ and $gg\m=t(g)\in \G_0$; see details in Subsection \ref{sub-grupoide}. If $\G_0$ has a unique element then $\G$ is a group. This algebraic version of groupoids motivated the authors of \cite{BP} to consider partial actions of groupoids on rings. In particular, it was defined in \cite{BP} the notion of Galois extension for partial actions of groupoids. A version of the Galois correspondence for global actions of groupoids on commutative rings was given in \cite{PT}.\vspace{0.07cm}

An special class of partial actions of connected groupoids was studied in \cite{BPP2}. This class was named {\it group-type} partial groupoid actions and this name is due to the fact that the partial skew groupoid ring associated can be realized as a partial skew group ring; see details in Theorem 4.4 of \cite{BPP2}. In this work we shall extend this notion for groupoids not necessarily connected. It is easy to construct examples of group-type partial actions of groupoids using the formulas given in (4) and (5) of \cite{BPP1}. In particular, every global groupoid action is a group-type partial action. \vspace{0.07cm}

The main contribution of this work is to give a Galois correspondence for group-type partial actions of groupoids. Precisely, let $\alpha=(S_g,\alpha_g)_{g\in \G}$ be a unital group-type partial action of a connected finite groupoid $\G$ on a comutative ring $S=\oplus_{y\in {\G}_0}S_y$. For each subgroupoid $\cH$ of $\G$, we consider $\af_{\cH}=(S_h,\alpha_h)_{h\in \cH}$ the partial action of $\cH$ on $S_{\cH}=\oplus_{y\in {\cH}_0}S_y$. Denote by $S^{\af_{\cH}}$ the subring of invariant elements (see Section 3). On the other hand, $\G_{T}$ denotes the set of elements of $\G$ that fix $T$, where $T$  is a subring of $S$ (see Section 4). Consider the set $\wSubg{\G}$ whose elements are wide subgroupoids $\cH$ of $\G$ such that $\alpha_{\cH}$ is group-type. Also, let $\mathfrak{B}(S)$ be the set of all subrings $T$ of $S$ which are $S^{\af_{\G}}$-separable, $\af$-strong and such that $\G_T=\cH$, for some $\cH\in \wSubg{\G}$. With this notation, we have the following Galois correspondence.\vspace{.15cm}

\noindent{\bf Theorem.}  {\rm \textbf{(Galois Correspondence)}} {\it Let $S$ be an $\alpha_{\G}$-partial Galois extension of $S^{\af_{\G}}$. There exists a bijective correspondence  between $\wSubg{\G}$ and  $\mathfrak{B}(S)$ given by ${\cH}\mapsto S^{\af_{\cH}}$ whose inverse is given by $T\mapsto \G_T$.}
\vspace{.15cm}

\noindent The Galois correspondence for not-necessarily connected groupoids follows from connected case; see Remark \ref{gen-case}. The previous theorem recover the Galois correspondence for global groupoid actions given in Theorem 4.6 (i) of \cite{PT}. \vspace{0.07cm} 

The organization of the paper is the following. As preliminaries, we recall the basic facts on groupoids, partial actions of groupoids on rings and Galois theory. In Section 3, given a unital partial action $\alpha=(S_g,\af_g)_{g\in \G}$  of a connected groupoid $\G$ on the ring $S=\oplus_{y\in\G_0}S_y$ and a subgroupoid ${\cH}$ of $\G$, we associate the subring $S^{\af_{\cH}}$ of invariant elements (by the partial action of ${\cH}$) of $S$. If the partial action $\af_{\cH}$ of $\cH$ on $S_{\cH}=\oplus_{y\in\cH_0} S_y$ (obtained by restriction of $\af$) is group-type then $S^{\af_{\cH}}$ is determined from invariant elements by the partial action of certain isotropy groups; see  Corollary \ref{prop-invariants-first}. In Section 4, given a subring $T$ of $S$, we consider the subset $G_T$ of the elements of $\G$ that acts trivially on $T$. The Proposition \ref{prop-importante} relates $G_T$ and $G(y)_{T_y}$, where $T_y=T1_y$ and $\G(y)$ is the isotropy group associated to $y\in \G_0$.
Sections 5 is dedicated to prove the theorem of Galois correspondence.     

\subsection*{Conventions}\label{subsec:conv}
Throughout this work, the rings are associative and not necessarily unital. An ideal of a ring means a two-sided ideal.
  
\section{Preliminaries}
The background about groupoids, partial groupoid actions and Galois theory that will be used in this work are presented here.

\subsection{Groupoids}\label{sub-grupoide}
We recall that a {\it groupoid} $\G$ is a small category in which every morphism is an isomorphism. Given $g\in \G$, the {\it source} and the {\it  target }  of $g$ will be denoted by $s(g)$ and $t(g)$, respectively. Also, $\G_0$ will denote the set of the objects of $\G$. Hence, we have maps $s,t:\G\to\G_0$.  Given $x\in \G_0$, $\id_x$ will denote the identity morphism of $x$.  Observe that $\id:\G_0\to \G$, given by $\id(x)=\id_x$, is an injective map and whence we identify $\G_0\subset \G$.
\vspace{0.05cm}

Given $x,y\in\G_0$ we set $\G(x,y):=\{g\in \G: s(g)=x\,\,\text{and}\,\,t(g)=y\}$. Notice that $\G(x):=\G(x,x)$ is a group which will be called the \emph{isotropy group associated to $x$.}
If $g\in \G(x,y)$ and $h\in \G(y,z)$ then $hg\in \G(x,z)$ will denote the composition of the elements $h$ and $g$. Also $\G^{(2)}=\{(g,h)\in\G\times\G \ :\ s(g)=t(h)\}$ denote the set of ordered pairs of $\G$ which are composable. It is clear that 
\begin{align}\label{useful-identities	}
& s(g)=g^{-1}g,& &t(g)=gg^{-1},& &s(gh)=s(h),& &t(gh)=t(g),&
\end{align}
for all $x\in \G_0$, $g\in \G$ and $(g,h)\in\G^{(2)}$.
\vspace{0.05cm}

A groupoid $\G$ is \emph{connected} if $\mathcal{\G}(x,y)\neq \emptyset$ for all $x,y\in \G_0$. It is well-known that any groupoid is a disjoint union of connected subgroupoids. Indeed, the equivalence relation on $\G_0$ given by 
\begin{align}\label{connected-equiv}
x\sim y \,\,\text{ if and only if }\, \G(x,y)\neq \emptyset, \qquad x,y\in \G_0,
\end{align}
induces the decomposition $\G=\dot\cup_{Y\in \G_0/\!\sim}\G_Y$ of $\G$ in connected components. For each equivalence class $Y\in\G_0/\!\!\sim$ we associated the connected component  $\G_Y$ of $\G$, which is a full connected subgroupoid of $\G$, in the following way.  The set of the objects of $\G_Y$ is $Y$ and $\G_Y(x,y)=\G(x,y)$ for all $x,y\in Y$. 
\vspace{0.05cm}

A connected groupoid is closely related to any of its isotropy group as we will see below. A proof of the next result can be seen, for instance, in $\S\,2$ of \cite{BPP1}. 

\begin{prop1}\label{group:connec}
	Let $\G$ be a connected groupoid, $x\in \G_0$, $\tau_y\in \G(x,y)$ for all $y\in \G_0$, $y\neq x$, $\tau_x=x$ and 
\begin{align}\label{group-element}
\tau(g):=\ta^{-1}_{t(g)}g\ta_{s(g)}\in \G(x),\quad g\in \G.
\end{align}	
Then $\psi:\G\to \G_0^2\times \G(x)$ given by $\psi(g)=((s(g),t(g)), \tau(g))$, for all $g\in \G$, is an isomorphism of groupoids.
\end{prop1}

\begin{remark}\label{brandt-groupoid}{\rm
The groupoid structure of $\G_0^2=\G_0\times \G_0$ is the usual, that is, $s(x,y)=x$ and  $t(x,y)=y$, for all  $x,y\in G_0$. The composition in $\G_0^2$ is given by $(y,z)(x,y)=(x,z)$ and the inverse of $(x,y)$ is $(y,x)$, for all $x,y,z\in \G_0$.}
\end{remark}

\subsection{Partial actions of groupoids}
We recall the basic notions an results related to partial actions of groupoids; the reference that will be used here is \cite{BP}.

\begin{def1}\label{def-partial-action}{\rm
A partial action of a groupoid $\G$ on  a ring $S$ is a family of pairs
$\af=(S_g,\af_g)_{g\in \G}$ that satisfies:
\begin{enumerate}[\rm (i)]
	\item for each $g\in \G$,  $S_{t(g)}$ is an ideal of $S$, $S_g$ is an ideal of $S_{t(g)}$ and $\af_g:S_{g\m}\to S_g$ is an isomorphism of rings,\vspace{0.05cm}
	
	\item for each $x\in \G_0$, $\af_x=\id_{S_x}$ is the identity map  of $S_x$, \vspace{0.05cm}
	\item $\af_h^{-1}(S_{g^{-1}}\cap S_h)\subseteq S_{(gh)^{-1}}$, for all $(g,h)\in \G^{(2)}$,\vspace{0.05cm}
	\item $\af_g(\af_h(a))=\af_{gh}(a)$, for all $a\in \af_h^{-1}(S_{g^{-1}}\cap S_h)$ and $(g,h)\in \G^{(2)}$.
\end{enumerate}}
\end{def1}

\noindent The items (iii) and (iv) above imply that $\af_{gh}$ is an extension of $\af_g\af_h$, that is,  $\alpha_{g}\alpha_{h}\leq \alpha_{gh}$, for all  $(g,h)\in \G^{(2)}$.\vspace{0.05cm}

\begin{def1}\label{def-global-action}{\rm
Consider a partial action $\alpha=(S_g,\af_g)_{g\in \G}$ of a groupoid $\G$ on  a ring $S$.
\begin{enumerate}[\rm (i) ]
	\item $\alpha $ is said {\it unital} if $S_g$ is a unital ring, for all $g\in \G$. In this case, for every $g\in \G$ there exists a central idempotent $1_g$ of $S$ such that $S_g=S1_g$,
	\item $\alpha$ is called {\it global} if $\af_g\af_h=\af_{gh}$, for all $(g,h)\in \G^{(2)}$.
\end{enumerate}} 
\end{def1}

\begin{remark}\label{obs-actions}
{\rm  Let $\af=(S_g,\af_g)_{g\in \G}$ be a partial action of a groupoid $\G$ on a ring $S$, $x\in \G_0$ and $\cH$ a subgroupoid of $\G$. Then:
\begin{enumerate}[\rm (i) ]
	\item the isotropy group $\G(x)$  acts partially on $S_x$ via $\af_{\G(x)}:=(S_g,\af_g)_{g\in \G(x)},$
	\item if $x\in {\cH}_0$ then  $\af_{{\cH}(x)}:=(S_h,\af_h)_{h\in {\cH}(x)}$ is a partial action of ${\cH}(x)$ on $S_x$; in this case $\af_{{\cH}(x)}$ is the restriction of $\alpha_{\G(x)}$ to ${\cH}(x)$,\vspace{0.05cm}
	\item  if $S=\bigoplus_{y\in \G_0}S_y$ and  ${\cH}$ is a  subgroupoid of $\G$ then we can consider the restriction $\alpha_{{\cH}}:=(S_h,\af_h)_{h\in {\cH}}$ of $\af$ to ${\cH}$ which is a partial action of ${\cH}$ on $S_H := \bigoplus_{z\in \cH_0}S_z$.
\end{enumerate}}
\end{remark}
Examples of partial actions of groupoids can be constructed via restriction of global actions. Indeed, let $\beta=(T_g,\bt_g)_{g\in \G}$ be a global action of a groupoid $\G$ on a ring $T$ and $S\subset T$ an ideal of $T$. For each $x\in \G_0$ and $g\in \G$, consider 
\begin{align*}
&S_x=S\cap T_x,& &S_g=S_{t(g)}\cap \beta_g(S_{s(g)}),& &\af_g=\beta_g|_{S_{g^{-1}}}.&
\end{align*}
Then $\af=(S_g,\af_g)_{g\in \G}$ is a partial action of $\G$ on $S$. \vspace{0.05cm}

Let $\alpha=(S_g,\alpha_g)_{g\in \G}$ be a partial action of a groupoid $\G$ on a ring $S$. It follows from Lemma 1.1 of \cite{BP} that 
\begin{align}\label{two-proper}	
&\alpha\m _g=\alpha_{g\m},\,\,\,\,g\in \G,& &\af_g(S_{g\m}\cap S_h)=S_g\cap S_{gh},\,\,\,\,(g,h)\in \G^{(2)}.&
\end{align}
Suppose that $\alpha$ is unital and assume that $S_g=S1_g$ where $1_g$ is a central idempotent of $S$, for all $g\in \G$. The equality from right side of \eqref{two-proper} implies that 
\begin{align}\label{eq-identities}
\af_g(\af_h(a1_{h\m})1_{g\m})=\af_{gh}(a1_{(gh)\m})1_g,\quad \text{ for all}\,\, (g,h)\in \G^{(2)},\,\,a\in S.
\end{align}

We end this section by proving a result that will be useful in what follows in this paper.

\begin{prop1}\label{prop-grouptype} Let $\af=(S_g,\af_g)_{g\in \G}$ be a partial action of a connected groupoid $\G$ on a ring $S$ and fix an element $\ta_{y,z}\in \G(y,z)$, for all $y,z\in \G_0$. If $S_{\ta^{-1}_{y,z}}=S_y$ and $S_{\ta_{y,z}}=S_z$ then 
	\begin{align}\label{ide-importante}
		S_g=S_{g\tau_{y,z}}, \,\text{ for all } g\in \G \text{ such that } s(g)=z.
	\end{align}
	In particular, we have that
	\begin{align}\label{eq:idempot}
		1_{g}=1_{g\tau_{y,z}}, \,\text{ for all } g\in \G \text{ such that } s(g)=z.
	\end{align}
\end{prop1}
\begin{proof}
	Since $\af$ is a partial action we have that $\af_{g\ta_{y,z}}$ is an extension of $\af_{g}\af_{\ta_{y,z}}$. Thus, the image $\operatorname{Im}(\af_{g}\af_{\ta_{y,z}})$ of $\af_{g}\af_{\ta_{y,z}}$ is a subset of $\operatorname{Im}(\af_{g\ta_{y,z}})=S_{g\ta_{y,z}}$. Since
	\begin{align*}
	\operatorname{Im}(\af_{g}\af_{\ta_{y,z}})&=\af_{g}(S_{g\m}\cap S_{\ta_{y,z}})=\af_{g}(S_{g\m}\cap S_z)\\
	&=\af_{g}(S_{g\m}\cap S_{s(g)})=\af_{g}(S_{g\m})=S_{g},
	\end{align*}
it follows that $S_g\subset S_{g\ta_{y,z}}$. Similarly, using that $\af_{g}=\af_{g\ta_{y,z}\ta\m_{y,z}}$ is an extension of $\af_{g\ta_{y,z}}\af_{\ta\m_{y,z}}$  and that $S_{\ta^{-1}_{y,z}}=S_y$, we obtain $S_{g\ta_{y,z}}\subset S_{g}$. 
\end{proof}

\subsection{Group-type partial actions}\label{subsec-group-type}  Given a connected groupoid $\G$ and $x\in \G_0$, consider the equivalence relation $\sim_x$ on $\cup_{y\in \G_0} \G(x,y)$ defined by
$$g\sim_x l\, \text{ if and only if there exists } z\in \G_0 \text{ such that }  g,l\in \G(x,z).$$
A transversal $\tau=\{\tau_y\}_{y\in \G_0}$ for $\sim_x$ such that $\tau_x=x$ will be called a \emph{transversal in $\G$ for $x$}. Hence,
$\tau_y:x\to y$ is a chosen morphism of $\G$, for each $y\in\G_0$. Notice that each transversal $\tau(x)$ in $\G$ for $x$ induces a groupoid epimorphism
$\tau:\G\to \G(x)$ given by
\begin{align}\label{ele-gr}
\tau(g)=\tau^{-1}_{t(g)}g\tau_{s(g)}, \quad \text{for all } g\in \G.
\end{align}
Notice that $\tau(\tau_y)=x$, for all $y\in \G_0$, and $\tau(l)=l$, for all $l\in \G(x)$.
\vspace{0.1cm}


The next definition was given in $\S\,3.2$ of \cite{BPP2} and  it will be fundamental for our purposes in this paper.

\begin{def1}\label{def-group-type}
	{\rm A partial action $\af=(S_g,\af_g)_{g\in \G}$ of a connected groupoid $\G$ on $S$ is called {\it  group-type} if there exist and element $x\in \G_0$ and a transversal $\tau=\{\tau_{y}\}_{y\in \G_0}$ in $\G$ for $x$ such that
		\begin{align}
			\label{cond1} S_{\tau\m_y}=S_x \ \ \text{and} \ \ S_{\tau_y}=S_y, \ \ \text{ for all } \ y\in\G_0.
	\end{align}}
\end{def1}

Any global groupoid action $\af=(S_g,\af_g)_{g\in \G}$ of $G$ on $S$ is group-type. In fact, Lemma 1.1 of \cite{BP} implies that $S_{g}=S_{t(g)}$, for all $g\in \G$. Hence \eqref{cond1} is satisfied and $\af$ is group-type. For examples of group-type partial actions that are not global we refer $\S\,3.2$ of \cite{BPP2}.

\begin{remark}\label{not-depend-x}
	{\rm Assume that $\af$ is a group-type partial action of $\G$ on $S$ and consider  $\tau=\{\tau_{y}\}_{y\in \G_0}$ a trans\-versal in $\G$ for $x$ such that \eqref{cond1} is true. Given $z\in \G_0$, fix the set $\gamma=\{\tau_y\ta\m_z\}_{y\in \G_0}$. By Remark 3.4 of \cite{BPP2}, $\gamma$ is a transversal in $\G$ for $z$ that also satisfies \eqref{cond1}. Thus, the notion of group-type partial action does not depend on the choice of the object $x$.}
	
\end{remark}

Now we will extend the notion of group-type partial actions for groupoids not necessarily connected. Let $\moL$ be a groupoid and $\moL=\cup_{Y\in \moL_0/\!\sim}\moL_Y$ its decomposition in disjoint connected components. Assume that $\af=(S_l,\af_l)_{l\in \moL}$ is a partial action of $\moL$ on a ring $S$ and that $S=\oplus_{z\in \moL_0}S_z$. In this case, for each $Y\in \moL_0/\!\!\sim$ we have a partial action $\af_Y$ of $\moL$ on $S_Y$, where
\begin{align} \label{partial-action-by-restriction-cc}
	&S_{Y}:=\mathsmaller{\bigoplus_{y\in Y}S_y},&&\quad\quad\af_Y:=(S_l,\af_l)_{l\in \moL_Y}.&
\end{align}

\begin{def1}\label{def-groupoid-tipogrupo}
	{\rm Let $\moL$ be a groupoid and $\af=(S_l,\af_l)_{l\in \moL}$ a partial action of $\moL$ on the ring $S=\oplus_{z\in \moL_0}S_z$.  We say that $\af$ is {\it group-type} if the partial action $\af_Y$ of the connected groupoid $\moL_Y$ on $S_Y$ is group-type (in the sense of Definition \ref{def-group-type}), for all  $Y\in \moL_0/\!\!\sim$.}
\end{def1}
The above definition will be used in Remark \ref{gen-case} to determine the Galois correspondence for any finite groupoid.

\subsection{Galois extension}\label{red-con} In what follows in this subsection, $\G$ denotes a finite groupoid and  $\alpha=(S_g,\alpha_g)_{g\in \G}$ is a unital partial action of $\G$ on a commutative ring $S$. We will assume $S_g=S1_g$, where $1_g$ is a central idempotent of $S$, for all $g\in \G$.
\begin{def1}\label{def-Galois-groupoide}
	{\rm
		The ring extension $R\subset S$ is called an  $\alpha_{\G}$-partial Galois extension if $R=S^{\af_{\G}}:=\{s\in S\,:\,\af_g(s1_{g\m})=s1_g,\text{ for all } g\in \G\}$ and there exist a positive integer $m$ and  elements $a_i,b_i\in S$, $1\leq i\leq m$, such that 
		\begin{align}\label{def-Galois}
			&\sum_{1\leq i\leq	m}a_i\alpha_g(b_i1_{g^{-1}})=\delta_{z,g}1_z,&  &\text{for all } z\in \G_0, \,\,g\in \G.& 
		\end{align}
		The set $\{a_i,b_i\}_{1\leq i\leq m}$ is called a  partial Galois coordinate system of $S$ over $R$. }
	
\end{def1}

Consider the decomposition $\G=\dot\cup_{Z_j\in \G_0/\!\sim}\G_{Z_j}$ of $\G$ in connected components, where $\G_0=Z_1\,\dot\cup\ldots \dot\cup\, Z_r$. Assume that  $S=\oplus_{z\in \G_0}S_z$ and  fix the notations: $\G_j:=\G_{Z_j}$, $S_j:=\oplus_{y\in Z_j}S_y$ and 
$\af_j$ the partial action of $\G_{j}$ on $S_j$ obtained by the restriction of $\af$ to $\G_{j}$, that is, $\af_{j}=(S_g,\af_{g})_{g\in \G_{j}}$. Moreover, for $g\in \G_{j}$, we will denote the partial bijection $(\af_j)_{g}:S_{g\m}\to S_g$ by  $\af_{j,g}$. \vspace{0.05cm}

One of the main purpose of this work is to present a Galois correspondence for partial actions of groupoids. The next result reduces this problem for connected groupoids; see Remark \ref{gen-case}.

\begin{prop1}\label{prop-reduction}
	Let $R=S^{\af_{\G}}$ and $R_j=S_j^{\af_j}$, for each $1\leq j\leq r$. Then $R \subset S$ is an $\af_{\G}$-partial Galois extension if and only if $R_j\subset S_j$ is an $\af_{j}$-partial Galois extension, for all $1\leq j\leq r$.
\end{prop1}
\begin{proof} Let $\{a_i, b_i\}_{1\leq i\leq m}$ be a partial Galois coordinate system of $S$ over $S^{\alpha}$. For each $Z_j \in G_0/\!\!\sim$ and $1\leq i\leq m$, we define $\tilde a_i = a_i1_{S_{j}}$ and $\tilde b_i = b_i1_{S_{j}}$. Notice that for every $h \in \G_{j}$,
	\begin{align*}
		\sum_{i=1}^m \tilde a_i \alpha_{j,h} (\tilde b_i1_{h^{-1}}) &  =  \sum_{i=1}^m  a_i1_{S_{j}} \alpha_{j,h}(b_i1_{S_{j}}1_{h^{-1}})= \sum_{i=1}^m  a_i \alpha_{j,h} (b_i1_{h^{-1}}) \\
		&  =  \sum_{i=1}^m  a_i \alpha_{h} (b_i1_{h^{-1}})= \delta_{y,h}1_y,
	\end{align*}
	for all $y \in Z_j$. Thus, $\{\tilde a_i,  \tilde b_i\}_{1\leq i\leq m}$ is a partial Galois coordinate system of $S_j$ over $S_j^{\alpha_j}$.	Conversely, let $\{a_{j,i}, b_{j,i}\}_{1\leq i\leq m_j}$ be a partial Galois coordinate system of $S_{j}$ over $S_j^{\alpha_{j}}$ for each $Z_j \in G_0/{\sim}$. Consider $m=\max\{m_j\,:\,1\leq j\leq r\}$ and fix $a_{j,i}:=0$ for all $m_j\leq i\leq m$. We also consider
	\[a_i:= \sum_{j=1}^r  a_{j,i},\qquad b_i:=\sum_{j=1}^r  b_{j,i},\qquad 1\leq i\leq m.\] 
	For each $g \in \G$, there exists a unique $k$ such that $g\in \G_{k}$. Then
	\begin{align*}
		\sum_{i=1}^m  a_i \alpha_{g}(b_i1_{g^{-1}}) & = \sum_{i=1}^m \sum_{j=1}^r \sum_{l=1}^r a_{j,i} \alpha_{g}(b_{l,i}1_{g^{-1}})= \sum_{i=1}^m   a_{k,i}\alpha_{g}( b_{k,i} 1_{g^{-1}}) \\
		&  = \sum_{i=1}^{m_k}   a_{k,i}\alpha_{g}( b_{k,i} 1_{g^{-1}})= \delta_{x,g}1_x,
	\end{align*}
	for all $x\in Z_k$. Since $\delta_{y,g}=0$ for all $y\in \G_0\setminus Z_k$, we conclude that $\{a_i, b_i\}_{1\leq i\leq m}$ is a partial Galois coordinate system of $S$ over $S^{\alpha}$.	
\end{proof}

\section{Subring of invariants}\label{invariants} 
In what follows in this paper, unless otherwise stated, $\alpha=(S_g,\alpha_g)_{g\in \G}$ is a unital group-type partial action of a finite connected groupoid $\G$ on a ring $S$. We will assume that  
\[S=\oplus_{y\in {\G}_0}S_y, \,\text{ where }\, S_g=S1_g. \]
Also, we will take a transversal  $\tau=\{\tau_{y}\}_{y\in \G_0}$ in $\G$ for a fixed object $x\in \G_0$.\vspace{0.05cm}

We recall from \cite{BP} that  an element $a\in S$ is called {\it $\af$-invariant} if $\af_g(a1_{g\m})=a1_g$, for all $g\in \G$. As in the previous section, the set of all $\af$-invariant elements of $S$ is denoted by 
\begin{align*}
	S^{\af_{\G}}:=\{a\in S\,:\,\af_g(a1_{g\m})=a1_{g},\,\text{ for all }\,g\in \G\}
\end{align*}
and it is a subring of $S$ which is called {\it subring of invariant elements of $S$}. 
Given $a=\sum_{y\in \G_0}a_y\in S$, with $a_y\in S_y$ for all $y\in \G_0$, there is a unique element $a_{x,y}\in S_x$ such that $a_y=\af_{\ta_y}(a_{x,y})$. Consequently, $a$ admits a unique decomposition as
\begin{align}
	\label{decomp-a} a=\sum_{y \in {\G}_0}\af_{\tau_y}(a_{x,y}).
\end{align}
In order to prove the next result we recall from \eqref{group-element} that
$\tau(g)=\ta^{-1}_{t(g)}g\ta_{s(g)}\in \G(x)$, for all  $g\in \G$.

\begin{prop1}\label{prop-invariants-first-one} 
Let $a\in S$ with decomposition as in \eqref{decomp-a}. Then $a\in S^{\alpha_{\G}}$  if and only if
	\begin{align}\label{cond-inv-prop}
		\af_{\ta(g)}\big(a_{x,s(g)}1_{\ta(g)\m}\big)=a_{x,t(g)}1_{\tau(g)}, \quad\text{for all } g\in {\G}. 
	\end{align}
\end{prop1}
\begin{proof}
	Let $g\in {\G}$. Notice that 	
	\begin{align*}
		\af_g(a1_{g\m})&=\af_g(\af_{\tau_{s(g)}}(a_{x,s(g)})1_{g\m})\\[.2em]
		&=\af_{g\tau_{s(g)}}(a_{x,s(g)}1_{(g\tau_{s(g)})\m})1_g\\
		&\stackrel{\eqref{ele-gr}}{=}\af_{\tau_{t(g)}\tau(g)}(a_{x,s(g)}1_{(\tau_{t(g)}\tau(g))\m})1_{\tau_{t(g)}}1_g\\
		&\stackrel{\eqref{eq-identities}}{=}\af_{\tau_{t(g)}}(\af_{\tau(g)}(a_{x,s(g)}1_{\tau(g)\m})1_x)1_g.
	\end{align*}
	Thus, $ a\in S^{\af_{\G}}$ if and only if    \[\af_{\tau_{t(g)}}(\af_{\tau(g)}(a_{x,s(g)}1_{\tau(g)\m})1_x)1_g =\af_g(a1_{g\m})=a1_g=\af_{\tau_{t(g)}}(a_{x,t(g)})1_g.\]
	From \eqref{eq:idempot} we have that $\af_{\tau\m_{t(g)}}(1_g)=1_{\tau(g)}$. Consequently $a\in S^{\af_{\G}}$ if and only if  
	$\af_{\tau(g)}(a_{x,s(g)}1_{\tau(g)\m})=a_{x,t(g)}1_{\tau(g)}$, for all $g\in {\G}$.
\end{proof}	
\vspace{0.05cm}

Using the previous result we have the following consequence.

\begin{cor1}\label{cor-invariants}
The map $\Phi_{\ta}:S_x^{\af_{\G(x)}}\to S^{\af_{\G}}$ given by 
	\begin{align}\label{eq:invariant}
		&\Phi_{\tau}(a)=\sum_{y\in {\G}_0} \af_{\tau_y}(a), \quad \text{ for all } a\in S_x^{\af_{\G(x)}},&
	\end{align}
is a ring isomorphism.
\end{cor1}
\begin{proof}
	Let $b=\sum_{y\in {\G}_0}\af_{\tau_{y}}(a)$ with $a\in S_x^{\af_{\G(x)}}$.  Then, for all $g\in {\G}$, we have
	\begin{align*}
		\af_g(b1_{g\m})&=\af_g(\af_{\tau_{s(g)}}(a)1_{g\m})\stackrel{\eqref{eq-identities}}{=}\af_{g\tau_{s(g)}}(a1_{(g\tau_{s(g)})\m})1_g\\
		&\stackrel{\eqref{ele-gr}}{=}\af_{\tau_{t(g)}\tau(g)}(a1_{(\tau_{t(g)}\tau(g))\m})1_g=\af_{\tau_{t(g)}\tau(g)}(a1_{(\tau_{t(g)}\tau(g))\m})1_{\tau_{t(g)}} 1_g\\
		&\stackrel{\eqref{eq-identities}}{=}\af_{\tau_{t(g)}}(\af_{\tau(g)}(a1_{\tau(g)\m})1_x)1_g=\af_{\tau_{t(g)}}(a1_{\tau(g)})1_g\\[.3em]
		&=\af_{\tau_{t(g)}}(a)\af_{\tau_{t(g)}}(1_{\ta(g)})1_g\stackrel{\eqref{eq:idempot}}{=}\af_{\tau_{t(g)}}(a)1_g=b1_g.
	\end{align*}
Hence, $\Phi_{\ta}(S_x^{\af_{\G(x)}})\subset S^{\af_{\G}}$.
	Conversely, let $a\in S^{\af_{\G}}$ with decomposition as in \eqref{decomp-a}. Since $\ta(\ta_{y})=x$, it follows from \eqref{cond-inv-prop} that  $a_{x,x}=a_{x,y}$, for all $y\in {\G}_0$. It is clear that $a_{x,x}\in S_x^{\af_{\G(x)}}$ which implies that $\Phi_{\ta}$ is surjective. Since $\Phi_{\ta}$ is an injective ring homomorphism, the result follows.
\end{proof}

\begin{remark}\label{not-depend-choose-object} {\rm Let $y\in \G_0$. By Remark \ref{not-depend-x} (i), $\gamma=\{\gamma_z\}_{z\in \G_0}$ with $\gamma_{z}=\ta_z\ta\m_y$, is a transversal for $y$ that satisfies \eqref{cond1}.  It follows that $S_y^{\af_{\G(y)}}\simeq S^{\af_{\G}}\simeq S_x^{\af_{\G(x)}}$, by Corollary \ref{cor-invariants}. Thus, 
$S_x^{\af_{\G(x)}}\simeq S_y^{\af_{\G(y)}}$, for all $y\in \G_0$.}		 
\end{remark}

Notice that Corollary \ref{cor-invariants} allows us to determine $S^{\af_G}$ from  $S_x^{\af_{G(x)}}$.  In order to use Corollary \ref{cor-invariants} for any subgroupoid (not-necessarily connected) of $\G$, we fix some extra notation. For a subgroupoid $\cH$ of $\G$, we put
\begin{align}\label{notation-rings}
&S_{{\cH}}:=\oplus_{z\in {\cH}_0}S_z,& &S^c_{{\cH}}:=S\setminus S_{\cH}=\oplus_{y\in {\cH}^c_0}S_y, \qquad \cH^c_0=\G_0\setminus\cH_0.&
\end{align}
According to Remark \ref{obs-actions} (iii), the restriction of $\af$ on ${\cH}$ gives a partial action $\af_{\cH}$ of ${\cH}$ on $S_{{\cH}}$, that is, $\af_{\cH}=(S_h,\alpha_h)_{h\in {\cH}}$.
\vspace{0.05cm}

\noindent{\bf Notation:} We will denote by $\Subg{\G}$ the set of all subgroupoids $\cH$ of $\G$ such that the partial action $\af_{\cH}$  of ${\cH}$ on $S_{{\cH}}$ is group-type (in the sense of Definition \ref{def-groupoid-tipogrupo}):
\[\Subg{\G}=\{\cH\,:\,\cH \text{ is a subgroupoid of }\G \text { and } \af_{\cH} \text{ is group-type}\}.\]

\vspace{0.05cm}

Consider a not-necessarily connected subgroupoid $\cH$  of $\G$ such that ${\cH}_0/\!\!\sim =\{Y_1,\ldots,Y_r\}$, where $\sim$ the equivalence relation on ${\cH}_0$ given in \eqref{connected-equiv}. Then
\begin{align} \label{subgroupoid-structure}
	&{\cH}_0=Y_1\,\dot{\cup}\ldots\dot{\cup}\, Y_r,& &{\cH}={\cH}_1\,\dot{\cup}\ldots\dot{\cup}\,{\cH}_r, \quad \text{with }\,{\cH}_j={\cH}_{Y_j}.&
\end{align}
Hence, for each $1\leq j\leq r$, we have the partial action $\af_{{\cH}_j}=(S_g,\alpha_g)_{g\in {\cH}_j}$ of ${\cH}_j$ on $S_j:=\oplus_{y\in Y_j} S_{y}$. 
Observe that, by Definition \ref{def-groupoid-tipogrupo}, $\cH\in\Subg{\G}$ if and only if $\af_{{\cH}_j}$ is a group-type partial action of ${\cH}_j$ on $S_j$, for all $j=1,\ldots,r$.

\begin{cor1}\label{prop-invariants-first}
Assume that $\cH\in \Subg{\G}$ has decomposition as in \eqref{subgroupoid-structure} and consider an object $y_j\in Y_j$, for all $j=1,\ldots,r$. Then
\begin{align}
S^{\af_{\cH}}=\big(\mathsmaller{ \bigoplus_{j=1}^r}\,S_j^{\af_{\cH_j}}\big)\mathsmaller{ \bigoplus} S^c_{{\cH}}  \simeq \big( \mathsmaller{ \bigoplus_{j=1}^r}\, S_{y_j}^{\af_{\cH_j(y_j)}}\big)\mathsmaller{ \bigoplus} S^c_{{\cH}} .	
\end{align}
\end{cor1}

\begin{proof}
Observe that $S=S_{1}\oplus\ldots\oplus S_{_r}\oplus  S^c_{{\cH}}$, where $S_{j}=\oplus_{z\in Y_j} S_z$. Then, an element $b\in S$ has a unique decomposition as $b=\sum_{j=1}^ra_j+a$, with $a_j\in S_j$ and $a\in  S^c_{{\cH}}$.
Clearly $b\in S^{\af_{\cH}}$ if and only if $a_j\in S_j^{\af_{\cH_j}}$. Consequently, $S^{\af_{\cH}}=\oplus_{j=1}^r S_j^{\af_{\cH_j}}\oplus S^c_{{\cH}}$ and the result follows from Corollary \ref{cor-invariants} and Remark \ref{not-depend-choose-object}.
\end{proof}

\begin{exe} \label{ex:invariant}
{\rm Let $\G=\{x, y, g,h,l,m,l^{-1},m^{-1}\}$ be the groupoid, where the set of objects is $\G_0=\{x,y\}$, and the composition rules are 
	\[g\in\G(x),\quad h\in\G(y),\quad l,m:x\to y\quad  g^2=x,\quad h^{2}=y,\quad lg=m=hl. \]
	Consider $S=\mathbb{C}e_1\oplus\mathbb{C}e_2\oplus\mathbb{C}e_3\oplus\mathbb{C}e_4$, where $\mathbb{C}$ denotes the complex number field, $e_ie_j=\delta_{i,j}e_i$ and $e_1+\ldots+e_4=1$. Define the following partial action $\alpha$ of $\G$ on $S$:
\begin{align*}
	S_x&=\mathbb{C}e_1\oplus\mathbb{C}e_2=S_{l\m},&  S_y&=\mathbb{C}e_3\oplus\mathbb{C}e_4=S_l, &\\[.2em]
	S_g&=\mathbb{C}e_1=S_{g\m}=S_{m\m},&  S_m&=S_h=\mathbb{C}e_3=S_{h\m}, &
	\end{align*}
	\begin{align*}
	\af_x&=\id_{S_x},& \af_g(ae_1)&=\overline{a}e_1,&  \af_m(ae_1)&=\overline{a}e_3, &\af_l(ae_1+be_2)&=ae_3+be_4, \\[.3em]
	\af_y&= \id_{S_y},& \af_h(ae_3)&=\overline{a}e_3,& \af_{m\m}(ae_3)&=\overline{a}e_1,&\af_{l\m}(ae_3+be_4)&=ae_1+be_2,
	\end{align*}
	where $\overline{a}$ denotes the complex conjugate of $a$, for every $a\in \mathbb{C}$. Observe that $\ta=\{\ta_x=x, \ta_y=l\}$ is a transversal in $\G$ for $x$  that  satisfies \eqref{cond1} and whence $\af$  is group-type. The set of subgroupoids of $\G$ are: 
	\begin{align*}
\cH_1&=\{x\}    	&	\cH_{4}&=\{x,g\},     	&   \cH_{7}&=\{x,y,h\},    &  \cH_{10}&=\{x,y,m,m\m\},\\
\cH_2&=\{y\}        &   \cH_{5}&=\{x,y,g\},     &   \cH_{8}&=\{x,y,g,h\},  &  \cH_{11}&=\G.\\         
\cH_{3}&=\{x,y\},   &   \cH_{6}&=\{y,h\},       &   \cH_{9}&=\{x,y,l,l\m\}, &          &           
	\end{align*}
Notice that  $\cH_j\in \Subg{\G}$ for all $j\ne 10$. We can determine the invariant subring $S^{\alpha_{\cH_j}}$ (for $j\neq 10$) of $S$ using the previous result. For instance, the decomposition of $\cH_{8}$ in disjoint connected components is $\cH_{8}=\cH_{4}\cup \cH_{6}$. Clearly, we have that  $S_{x}^{\cH_{4}(x)}=\mathbb{R}e_1\oplus \mathbb{C}e_2$ and  $S_{y}^{\cH_{6}(y)}=\mathbb{R}e_3\oplus \mathbb{C}e_4$. From Corollary \ref{prop-invariants-first} follows that $S^{\alpha_{H}}\simeq \mathbb{R}e_1\oplus \mathbb{C}e_2\oplus \mathbb{R}e_3\oplus\mathbb{C}e_4$. Note also that $S^{\cH_1}=S^{\cH_2}=S^{\cH_3}=S$. }
\end{exe}

\section{Subgroupoid arising from invariant subring}

In the previous section we associated for each subgroupoid ${\cH}$ of $\G$ the subring $S^{\alpha_{H}}$ of $S$. The purpose here is to present a correspondence in the reverse direction. 

From now on, we will assume that $1_g\neq 0$, for all $g\in \G$. By Example \ref{ex:invariant} it is clear that the correspondence $\cH\mapsto S^{\af_{\cH}}$ from  $\Subg{\G}$ to the set of subrings of $S$ is not injective. In order to obtain a one-to-one correspondence between subgroupoids in $\Subg{\G}$ and certain subalgebras of $S$,  we will introduce the following notation.

\smallbreak

\noindent{\bf Notation:} We set $$\wSubg{\G} := \{ \cH \in \Subg{\G} : \cH_0=\G_0\}.$$
the set of the wide subgroupoids $\cH$ contained in $\Subg{\G}$.  

\smallbreak

Throughout this section, $\cH\in \wSubg{\G}$  admits a decomposition as in \eqref{subgroupoid-structure}, that  is, ${\cH}={\cH}_1\,\dot{\cup}\ldots\dot{\cup}\,{\cH}_r$, where ${\cH}_j={\cH}_{Y_j}$ and ${\cH}_0=Y_1\,\dot{\cup}\ldots\dot{\cup}\, Y_r=\G_0$. Moreover, \linebreak $\af_{{\cH}_j}=(S_g,\alpha_g)_{g\in {\cH}_j}$ is the partial action of $\cH_j$ on $S_j:=\oplus_{y\in Y_j} S_{y}$ obtained by restriction. For each $1\leq j\leq r$, we take  $y_j\in Y_j$ and we consider a transversal $\tau_j=\{\ta_{j,z}:y_j\to z\}_{z\in Y_j}$ in $\cH_j$ for $y_j$  satisfying \eqref{cond1}.

Given a subring $T$ of $S$, we will denote by ${\G}_T$ the set of elements of $\G$ which act trivially on $T$, that is,
\begin{align*}
	{\G}_T:=\{g\in \G\,:\,\af_{g}(t1_{g\m})=t1_g,\,\text{ for all } t\in T\}.
\end{align*} 

\begin{remark}
	{\rm Notice that $\G_0\subset {\G}_T$ because $\af_y=\id_{S_{y}}$, for all $y\in \G_0$. In general, ${\G}_T$ is not a subgroupoid of $\G$. In fact, let $\G$, $S$ and $\alpha$ be as in Example \ref{ex:invariant}. Consider the subalgebra $T = \mathbb{R}(e_1 + e_3)\oplus \mathbb{C}e_2\oplus  \mathbb{C}e_4$ of $S$. It is immediate to verify that $m, g \in {\G}_T$ but $l\notin {\G}_T$. Since $mg = l$, it follows that ${\G}_T$ is not a subgroupoid of $\G$. However, if $\af$ is global then $\G_T$ is a subgroupoid of $\G$, by \cite[Lemma 2.1]{PT}}.
\end{remark}

\begin{lem1}\label{lem-invari-subgroupoid}
	Let $T=S^{\af_H}$, $T_{y_j}=S_{y_j}^{\af_{\cH_j(y_j)}}$ for all $j=1,\ldots,r$ and $g\in \G$. Then, $g\in \G_T$ if and only if  there exists $1\leq k\leq r$ such that $s(g),\,t(g)\in Y_k$ and $\ta_k(g)\in \G(y_k)_{T_{y_k}}$.
\end{lem1}

\begin{proof}
	By Corollary \ref{prop-invariants-first}, $T=T_1\oplus\ldots\oplus T_r$ with $T_j=S_j^{\af_{\cH_j}}$. Moreover, we have that $T_j=\big\{\sum_{y\in Y_j} \af_{\ta_{j,y}}(a_j)\,:\,a_j\in T_{y_j} \big\}$. Assume that $g\in \G_T$ and consider $1\leq k\leq r$ such that $s(g)\in Y_k$. Suppose that $t(g)\notin Y_k$. Since $g\in \G_T$ and $1_{k}=\sum_{y\in Y_k}1_y\in T$, it follows that $0=1_{k}1_{t(g)}1_g=1_{k}1_{g}=\af_g(1_{k}1_{g\m})=\af_g(1_{s(g)}1_{g\m})=1_g$, a contradiction. So $s(g),t(g)\in Y_k$, so that $\ta_k(g) = \ta\m_{k,t(g)}g\ta_{k,s(g)}\in \G(y_k)$. In order to prove that $\ta_k(g) \in \G(y_k)_{T_{y_k}}$ we take an arbitrary element $a_k\in T_{y_k}$ and consider $t=\sum_{y\in Y_k}\af_{\tau_{k,y}}(a_k)\in T_k$. Notice that 
		
	\begin{align*}
		\af_g(t1_{g\m})&=\af_{\ta_{k,t(g)}{\tau_k}(g)\ta\m_{k,s(g)}}(\af_{\ta_{k,s(g)}}(a_k)1_{g\m})\\
		&\stackrel{\eqref{eq-identities}}{=}\af_{\ta_{k,t(g)}{\tau_k}(g)}(a_k1_{(\ta_{k,t(g)}{\tau_k}(g))\m})1_{g}\\[.3em]
		&=\af_{\ta_{k,t(g)}{\tau_k}(g)}	(a_k1_{(\ta_{k,t(g)}{\tau_k}(g))\m})1_{\ta_{k,t(g)}}1_{g}\\
		&\stackrel{\eqref{eq-identities}}{=}\af_{\ta_{k,t(g)}}(\af_{{\tau_k}(g)}(a_k1_{{\tau_k}(g)\m})1_{y_k})1_{g}.
	\end{align*}
	Thus  $\af_{\ta_{k,t(g)}}(\af_{{\tau_k}(g)}(a_k1_{{\tau_k}(g)\m})1_{y_k})1_{g} = \af_g(t1_{g\m}) = t1_g = \af_{\ta_{k,t(g)}}(a_k)1_g$.  Applying now $\af\m_{\ta_{k,t(g)}}$ to the first and last terms of the previous identity we obtain
	\[a_k\af_{\tau\m_{k,t(g)}}(1_g)=\af_{{\tau_k}(g)}(a_k1_{{\tau_k}(g)\m})\af_{\tau\m_{k,t(g)}}(1_g).\]	
From \eqref{eq:idempot}, $\af_{\tau\m_{k,t(g)}}(1_g)=1_{{\tau_k}(g)}$ which implies that $\af_{{\tau_k}(g)}(a_k1_{{\tau_k}(g)\m})=a_k1_{{\tau_k}(g)}$. Hence ${\tau_k}(g)\in \G(y_k)_{T_{y_k}}$.   \vspace{0.05cm}
	
Conversely, consider $t\in T$. Hence $t=t_1+\ldots t_r$ with $t_j\in T_j$. Using again that $g=\ta_{k,t(g)}\ta_k(g)\ta\m_{k,s(g)}$ together with our hypotheses we obtain
	\begin{align*}
		\af_g(t1_{g\m})&=\af_g(t_k1_{g\m})=\af_{\ta_{k,t(g)}\ta_k(g)\ta\m_{k,s(g)}}(\af_{\ta_{k,s(g)}}(a_k)1_{g\m})\\
		&\stackrel{\eqref{eq-identities}}{=}\af_{\ta_{k,t(g)}\ta_k(g)}(a_k1_{(\ta_{k,t(g)}\ta_k(g))\m})1_g=\af_{\ta_{k,t(g)}\ta_k(g)}(a_k1_{(\ta_{k,t(g)}\ta_k(g))\m})1_{t(g)}1_g\\
		&\stackrel{\eqref{eq-identities}}{=}\af_{\ta_{k,t(g)}}(\af_{\ta_k(g)}(a_k1_{\ta_k(g)\m})1_{y_k})1_g=\af_{\ta_{k,t(g)}}(a_k1_{\ta_k(g)})1_g\\
		&=\af_{\ta_{k,t(g)}}(a_k)1_{\ta_{k,t(g)}\ta_k(g)}1_g\stackrel{(\ast)}{=}\af_{\ta_{k,t(g)}}(a_k)1_g=t_k1_g=t1_{g},
	\end{align*}
	where $(\ast)$ follows from \eqref{ide-importante} because $\ta_{k,t(g)}\ta_k(g)=g\ta_{k,s(g)}$. Therefore, $g\in \G_T$.	
\end{proof}


\begin{prop1}\label{prop-importante}
	Let  $T$ and $T_{y_j}$ as in the previous lemma. The following assertions are satisfied:
	\begin{enumerate}[\rm (i)]
		\item $\G_T$ is a wide subgroupoid of $\G$ if and only if $\G(y_j)_{T_{y_j}}$ is a subgroup of $\G(y_j)$, for all $1\leq j\leq r$,\vspace{0.05cm}
		\item $\G_T=\cH$ if and only if  $\G(y_j)_{T_{y_j}}=\cH_j(y_j)$, for all $1\leq j\leq r$. 
	\end{enumerate}		
\end{prop1}
\begin{proof}
	\noindent (i) Assume that $\G_T$ is a wide subgroupoid of $\G$. Let $1\leq j\leq r$ and $l_1,l_2\in \G(y_j)_{T_{y_j}}$. Since $\ta_{j}(l_1)=l_1$ and $\ta_{j}(l_2)=l_2$, it follows from Lemma \ref{lem-invari-subgroupoid} that $l_1,l_2\in \G_T$. Hence, $l_1l_2\in \G_T$. Using again Lemma \ref{lem-invari-subgroupoid} we obtain that $l_1l_2=\ta_{j}(l_1l_2)\in \G(y_j)_{T_{y_j}}$. Also, it is clear that $l\m_1\in \G(y_j)_{T_{y_j}}$. Conversely, let $g,h\in \G_{T}$ such that $s(g)=t(h)$. From Lemma \ref{lem-invari-subgroupoid}, there are $1\leq k_1,k_2\leq r$ such that $s(g), t(g)\in Y_{k_1}$, $s(h), t(h)\in Y_{k_2}$, $\ta_{k_1}(g)\in \G(y_{k_1})_{T_{y_{k_1}}}$ and $\ta_{k_2}(h)\in \G(y_{k_2})_{T_{y_{k_2}}}$. From $s(g)=t(h)$ follows that $k_1=k_2$ and cosequently $\ta_{k_1}(gh)=\ta_{k_1}(g)\ta_{k_1}(h)\in \G(y_{k_1})_{T_{y_{k_1}}}$, since $\G(y_{k_1})_{T_{y_{k_1}}}$ is a group by hipothesis. Then, Lemma \ref{lem-invari-subgroupoid} ensures that $gh\in \G_T$. It is trivial to check that $g\m \in \G_T$.  \vspace{0.05cm}
	
	\noindent (ii) Consider $1\leq j\leq r$ and assume that $\G_T=\cH$. Note that $\cH_j(y_j)\subset\G(y_j)_{T_{y_j}}$. Given $l\in\G(y_j)_{T_{y_j}}$, we have that $l=\ta_{j}(l) \in\G(y_j)_{T_{y_j}}$. By Lemma \ref{lem-invari-subgroupoid}, $l\in \G_T=\cH$. Hence, $l\in \cH(y_j)$. For the converse, notice that $\cH \subset \G_T$. If $g\in \G_T$ then by Lemma \ref{lem-invari-subgroupoid}, there exists $1\leq k\leq r$ such that  $s(g),\,t(g)\in Y_k$ and $\ta_k(g)\in \G(y_k)_{T_{y_k}}$. Thus $\G(y_k)_{T_{y_k}}=\cH_k(y_k)$ implies that $g=\ta_{k,t(g)}\ta_k(g)\ta\m_{k,s(g)}\in \cH_k\subset \cH$.	
\end{proof}


\begin{exe}\label{ex:groupoid-sec4}{\rm
		Let $\G=\{g,g^2,h,h^2,l,l\m,m,m\m,n,n\m\}$ be the groupoid with set of objects $\G_0=\{x,y\}$. The composition in $\G$ is given by
		\begin{align*}
			g&\in \G(x),\quad h\in \G(y),& &\,\,\,g^3=x,& &\,h^3=y,& \\
			&l,m,n\in\G(x,y),& lg&=m=hl,& mg&=n=hm.&
		\end{align*} 
		The groupoid $\G$ is illustrates in the following diagram 
		\[\xymatrix{& x\ar[d]^{n}\ar[r]^{g}  &x\ar[d]^{m} \ar[r]^{g} &x\ar[d]^{l}\\
			& y & y\ar[l]^{h}&\ar[l]^{h} y}\\ \]
		Let  $S:=\oplus_{i=1}^6k e_i$, where  $e_ie_j=\delta_{i,j}e_i$ and $\sum_{i=1}^{6}e_i=1$. Define the following partial action $\alpha$ of $\G$ on $S$:
		\begin{align*}
			&\af_x=\id_{S_x},\,\,\,S_{x}=k e_1\oplus&  &\!\!\!\!k e_2\oplus  k e_3,&  &\af_y=\id_{S_y},\,\,\,S_{y}=k e_4\oplus k e_5\oplus  k e_6, & \\
			&\af_l: k e_1\oplus k e_2\oplus k e_3\to& &\!\!\!\!k e_4\oplus k e_5\oplus k e_6,& &\af_l(e_i)=e_{i+3},\,\,\,i=1,2,3.& \\
			&\af_g: k e_1\to k e_2,& &\!\!\!\!\af_g(e_1)=e_2,& & \af_{g^2}=\af\m_g,& \\
			&\af_h: k e_4\to k e_5,& &\!\!\!\!\af_h(e_4)=e_5,& &\af_{h^2}=\af\m_h,&\\
			&\af_m: k e_1\to k e_5,& &\!\!\!\!\af_m(e_1)=e_5,& &\af_{m\m}=\af\m_m&\\
			&\af_n: k e_2\to k e_4,& &\!\!\!\!\af_n(e_2)=e_4,& &\af_{n\m}=\af\m_n.&
		\end{align*}
		Note that  $\af$ is a group-type partial action since $\tau=\{\tau_x=x,\tau_y=l\}$ is a transversal for $x$ that satisfies \eqref{cond1}. Take the subgroupoid $\cH=\{x,y,g,g^2,h,h^2\}$ of $\G$ and $T=S^{\af_{\cH}}$. Thus $\cH=\cH_1\dot{\cup} \cH_2$, where $\cH_1=\{x,g,g^2\}$ and $\cH_2=\{y,h,h^2\}$. Hence, it is clear that $\cH\in \wSubg{\G}$. Moreover, we have that $T=T_1\oplus T_2$ with $T_1=S_x^{\af_{\cH_{1}}}$ and $T_2=S_y^{\af_{\cH_{2}}}$. Notice that \[T_x=S_x^{\af_{\cH_{1}(x)}}=k(e_1+e_2)\oplus ke_3,\qquad T_y=S_y^{\af_{\cH_{2}(y)}}=k(e_4+e_5)\oplus ke_6.\]
		Thus, from Corollary \ref{cor-invariants} follows that $T=k(e_1+e_2)\oplus ke_3 \oplus  k(e_4+e_5)\oplus ke_6$. Since $\G(x)_{T_x}=\cH_1(x)=\cH_{1}$ and  $\G(y)_{T_y}=\cH_2(y)=\cH_{2}$, it follows from Proposition \ref{prop-importante}  that $\G_{T}=\cH$.	} 
\end{exe}

\section{Galois correspondence}
Throughout this section, $S=\oplus_{y\in {\G}_0}S_y$ is a commutative ring. As in the previous section, we shall assume that $\cH\in \wSubg{\G}$  admits a decomposition in connected components as  ${\cH}={\cH}_1\,\dot{\cup}\ldots\dot{\cup}\,{\cH}_r$, where ${\cH}_j={\cH}_{Y_j}$,  and ${\cH}_0=Y_1\,\dot{\cup}\ldots\dot{\cup}\, Y_r=\G_0$. Also, $\af_{{\cH}_j}=(S_g,\alpha_g)_{g\in {\cH}_j}$ is the partial action of $\cH_j$ on $S_j:=\oplus_{y\in Y_j} S_{y}$ and, for each $1\leq j\leq r$, we will take  $y_j\in Y_j$ and $\tau_j=\{\ta_{j,z}\,:\, z\in Y_j\}$ a transversal in $\cH_j$ for $y_j$  that satisfies \eqref{cond1}. \smallbreak

\subsection{The main theorem}
Our purpose in the sequel is to present a Galois correspondence between subgroupoids of $\cH$ and subrings of $S$. We start with the following.

\begin{lem1}\label{obs-Galois-compon}
	If $S^{\af_{\G}}\subset S$ is an $\af$-partial Galois extension, then $S_y^{\af_{\G(y)}}\subset S_y$ is an $\af_{\G(y)}$-partial Galois extension, for all $y\in\G_0$.
\end{lem1}
\begin{proof}
	Let $y\in \G_0$  and assume that $S$ is an $\af$-partial Galois extension of $R=S^{\af_{\G}}$. Then, there are $m$ and  elements $a_i,b_i\in S$, $1\leq i\leq m$, that satisfy \eqref{def-Galois}. Fix $a_{y,i}=a_i1_y$ and $b_{y,i}=b_i1_y$, for each $1\leq i\leq m$. It is immediate to check that $\{a_{y,i},b_{y,i}\}_{1\leq i\leq m}$ is a partial Galois coordinate system of $S_y$ over $R_y=S_y^{\G(y)}$.
\end{proof}

We recall that a  unital ring extension $R\subset T$ is called {\it separable} if the multiplication map $m:T\otimes_R T\to T$ is a splitting epimorphism of $T$-bimodules. This is equivalent to saying that there exists an element $e\in T\otimes_RT$ such that $te=et$, for all $t\in T$, and $m(e)=1_T$. Such an element $e$ is usually called {\it an idempotent of separability} of $T$ over $R$.

\begin{lem1}\label{lem-separavel}  
	Let $\cH\in \wSubg{\G}$, $T=S^{\af_H}$ and $R=S^{\af_G}$. For each $1\leq j\leq r$, consider $T_{y_j}=S_{y_j}^{\af_{\cH_j(y_j)}}$ and $R_{y_j}=S_{y_j}^{\af_{\G(y_j)}}$.
	Then the following statements are equivalent:
	\begin{enumerate}[\rm (i)] 
		\item $R\subset T$ is separable,\vspace*{.1cm}
		\item $R_{y_j}\subset T_{y_j}$ is separable, for all $1\leq j\leq r$.
	\end{enumerate}
\end{lem1}
\begin{proof}
	By Corrolary \ref{prop-invariants-first}, 
	$S^{\af_{\cH}}\simeq \bigoplus_{j=1}^r\, S_{y_j}^{\af_{\cH(y_j)}}$. Also, by Remark \ref{not-depend-choose-object} we have that $S^{\af_{\G}}\simeq R_{y_j}$, for all $j=1\ldots r$.
Hence the result is an immediate consequence from Proposition 1.7 (c) of \cite{K}.
\end{proof}

For the next definition $\G$ can be any groupoid (not necessarily connected). 

\begin{def1}\label{alpha-strong}{\rm
		Let $T$ be a subring of $S$ and $T_y:=T1_y$, for all $y\in \G_0$. 
		\begin{enumerate}[\rm (i)] 
			\item The subring $T_y$ of $S_y$ will be called {\it $\af_{\G(y,z)}$-strong} if for any $g,h\in \G(y,z)$ such that $g\m h\notin\G(y)_{T_y}$ and for any non-zero idempotent  $e\in S_g\cup S_h$, there exists $t_y\in T_y$ such that $\af_g(t_y1_{g\m})e\neq\af_h(t_y1_{h\m})$.\smallbreak
			\item  We shall say that $T$ is {\it $\af$-strong} if $T_y$ is $\af_{\G(y,z)}$-strong for all $y,z\in \G_0$.
	\end{enumerate}}
\end{def1}

The next result characterizes when $T=S^{\af_{\cH}}$ is $\af$-strong in the case where ${\cH}\in \wSubg{\G}$.

\begin{prop1}\label{prop-strong} Let ${\cH}\in \wSubg{\G}$ be a wide subgroupoid of $\G$, $T=S^{\af_{\cH}}$ and $T_{y_j}=S_{y_j}^{\af_{\cH_j(y_j)}}$ for all $1\leq j\leq r$. The following statements are equivalent:
	\begin{enumerate}[\rm (i)] 
		\item $T$ is $\af$-strong, \smallbreak
		\item for any $g,h\in \G$ such that $t(g)=t(h)$ and $g\m h\notin \G_T$ and for any non-zero idempotent  $e\in S_g\cup S_h$, there exists $t\in T$ such that $\af_g(t1_{g\m})e\neq \af_h(t1_{h\m})e$,\smallbreak 
		\item $T_{y_j}$ is $\af_{\G(y_j)}$-strong for all $1\leq j\leq r$.
	\end{enumerate}	
\end{prop1}

\begin{proof}
	
	Since $\af$ is group-type, there are $x\in \G_0$ and a transversal  $\ta=\{\ta_z\}_{z\in \G_0}$ in $\G$  for $x$  that satisfies \eqref{cond1}. Given $g\in \G$ with $s(g)\in Y_k$ and $t(g)\in Y_l$ we put 
	\[u(g):=\ta\m_{l,t(g)}g\ta_{k,s(g)}\in \G(y_k,y_l),\qquad v(g):=\ta_{y_k}\ta\m_{y_l}u(g)\in \G(y_k).\]
	
	\noindent (i)$\Rightarrow$(ii) Let $g,h\in \G$ with $t(g)= t(h)\in Y_l$, for some $1\leq l\leq r$, such that $g\m h\notin \G_T$. Arguing by contradiction, we will assume in addition that there exists a non-zero idempotent $e\in S_g\cup S_h$ of $S$ such that $\af_g(t1_{g\m})e=\af_h(t1_{h\m})e$, for all $t\in T$. 
	\smallbreak
	
	\smallbreak
	
	\noindent{\bf Claim 1.} There is $1\leq k\leq r$ such that $s(g),\,s(h)\in Y_k$.\vspace{.05cm}	
	
	\noindent Suppose that $s(g)\in Y_i$ and $s(h)\in Y_j$ with $1\leq i\neq j\leq r$. If $e\in S_g$ then we consider the following element $t_0=\sum_{y\in Y_i}1_{y}\in T$. Thus \[1_g=1_ge=\af_g(1_{g\m})e=\af_g(1_{s(g)}1_{g\m})e=\af_g(t_01_{g\m})e=\af_{h}(t_01_{h\m})e=0,\]
	which is an absurd. Similarly, if $e\in S_h$ then for $t_1=\sum_{y\in Y_j}1_{y}\in T$ we obtain that $1_h=1_he=\af_h(t_11_{h\m})e=\af_g(t_11_{g\m})e=0$. Hence, the claim follows.
	
	\smallbreak\smallbreak	
	
	\noindent{\bf Claim 2.} $\tau_k(g\m h)=u(g)\m u(h)\notin {\G(y_k)}_{T_{y_k}}$.
	
	\noindent This claim follows directly from Claim 1 and Lemma \ref{lem-invari-subgroupoid}.
	\smallbreak\smallbreak
	
	\noindent{\bf Claim 3.}  $f=\af_{\tau\m_{l,t(g)}}(e)$ is a non-zero idempotent of $S$ that belongs to $S_{u(g)}\cup S_{u(h)}$. \smallbreak
	
	\noindent Observe that $e\in S_{g}\cup S_h\subset S_{t(g)}=S_{t(h)}$. Suppose that $e\in S_g$. By \eqref{ide-importante} we have that  $S_g=S_{g\tau_{k,s(g)}}$ and consequently
	\begin{align*}
		f\in \af_{\tau\m_{l,t(g)}}(S_{g\tau_{k,s(g)}})=\af_{\tau\m_{l,t(g)}}(S_{g\tau_{k,s(g)}}\cap S_{t(g)})=S_{\tau\m_{l,t(g)}g\tau_{k,s(g)}}\cap S_{y_l}=S_{u(g)}.
	\end{align*}
	Similarly, if $e\in S_h$ then $f=\af_{\tau\m_{l,t(h)}}(e)\in S_{u(h)}$.
	\smallbreak\smallbreak
	
	\noindent{\bf Claim 4.} $\af_{u(g)}(a1_{u(g)\m})f=\af_{u(h)}(a1_{u(h)\m})f$, for all $a\in T_{y_k}$.   \smallbreak
	
	\noindent  Let $a\in T_{y_k}$. Then, by \eqref{eq:invariant}, $t:=\sum_{y\in Y_k}\af_{\tau_{k,y}}(a)\in S_k^{\af_{\cH_k}}$. Firstly, observe that
	\begin{align*}
		\af_g(t1_{g\m})e&=\af_g(\af_{\tau_{k,s(g)}}(a)1_{g\m})e\stackrel{\eqref{eq-identities}}{=}\af_{g\tau_{k,s(g)}}\big(a1_{(g\tau_{k,s(g)})\m}\big)1_ge\\
		&=\af_{\tau_{l,t(g)}u(g)}\big(a1_{(\tau_{l,t(g)}u(g))\m}\big)1_{t(g)}1_ge\\
		&\stackrel{\eqref{eq-identities}}{=}\af_{\tau_{l,t(g)}}\big(\af_{u(g)}(a1_{u(g)\m})\big)1_ge.
	\end{align*}
	Similarly, one can obtain that  $\af_h(t1_{h\m})e=\af_{\tau_{l,t(h)}}\big(\af_{u(h)}(a1_{u(h)\m})\big)1_he$.
By hypothesis, $\af_g(t1_{g\m})e=\af_h(t1_{h\m})e$. Thus, we have that
	\[\af_{\tau_{l,t(g)}}\big(\af_{u(g)}(a1_{u(g)\m})\big)1_ge=\af_{\tau_{l,t(h)}}\big(\af_{u(h)}(a1_{u(h)\m})\big)1_he.\]
Using now that $t(g)=t(h)$ we can apply $\af_{\tau\m_{l,t(g)}} = \af_{\tau\m_{l,t(h)}}$ in both sides of the last equality to obtain that 
	\begin{align*}
		\af_{u(g)}(a1_{u(g)\m})\af_{\tau\m_{l,t(g)}}(1_g)f 
		&=\af_{u(h)}(a1_{u(h)\m})\af_{\tau\m_{l,t(h)}}(1_h)f.
	\end{align*} 
	Finally, one can observe that $\af_{\tau_{l,t(g)}}(1_{u(g)})=1_{g\ta_{k,s(g)}}\stackrel{\eqref{ide-importante}}{=}1_g$. Hence $\af_{\tau\m_{l,t(g)}}(1_g)=1_{u(g)}$. Similarly, $\af_{\tau\m_{l,t(h)}}(1_h)=1_{u(h)}$ and the proof of Claim 4 is finished.
	\smallbreak
	
	From Claims $1$, $2$, $3$ and $4$ it follows that $T_{y_k}$ is not $\af_{\G(y_k,y_l)}$-strong and whence $T$ is not $\af$-strong, a contradiction. 
	
	\smallbreak

	\noindent (ii)$\Rightarrow$(iii) Consider $1\leq j\leq r$. Let  $g,h\in \G(y_j)$ such that $g\m h\notin \G(y_j)_{T_{y_j}}$ and $e\in S_{g}\cup S_{h}$ a non-zero idempotent of $S$. Clearly, $g\m h\notin\G_{T}$. Then, there exists $t\in T$ such that $\af_{g}(t1_{g\m})e\neq \af_h(t1_{h\m})e$. Take $t_{y_j}:=t1_{y_j}\in T_{y_j}$ and notice that $t1_{g\m}=t_{y_j}1_{g\m}$ and $t1_{h\m}=t_{y_j}1_{h\m}$. Thus, $\af_{g}(t_{y_j}1_{g\m})e\neq \af_h(t_{y_j}1_{h\m})e$ and whence $T_{y_j}$ is $\af_{\G(y_j)}$-strong.

	\smallbreak
	
	\noindent (iii)$\Rightarrow$(i) Let $y,z\in \G_0$, $T_y=T1_y$ and $T_z=T1_z$. Consider $g,h\in \G(y,z)$ such that $g\m h\notin\G(y)_{T_y}$ and  let $e\in S_g\cup S_h$ be a non-zero idempotent of $S$. Assume that $y\in Y_k$ and $z\in Y_l$. Then, as in Claim 3, $f:=\af_{\tau\m_{l,z}}(e)\in S_{u(g)}\cup S_{u(h)}$ is a non-zero idempotent of $S$. Hence, it is clear that  $f':=\af_{\ta_{y_k}\ta\m_{y_l}}(f)\in S_{v(g)}\cup S_{v(h)}$ is a nonzero-idempotent of $S$. By Lemma \ref{lem-invari-subgroupoid}, $\ta_{k}(g\m h)=u(g)\m u(h)\notin \G(y_k)_{T_{y_k}}$. Since $v(g)\m v(h)=u(g)\m u(h)$ it follows that $v(g)\m v(h)\notin \G(y_k)_{T_{y_k}}$. Therefore, there is $t_{y_k}\in T_{y_k}$ such that
	\begin{align}\label{eq-strong}
		\af_{v(g)}(t_{y_k}1_{v(g)\m})f'\neq \af_{v(h)}(t_{y_k}1_{v(h)\m})f'.	
	\end{align}
	It follows from Corollary \ref{cor-invariants} that $\af_{\ta_{k,y}}(T_{y_k})=T_y$ and then $t_y:=\af_{\ta_{k,y}}(t_{y_k})\in T_y$. Also, we have
	\begin{align*}
		\af_{v(g)}(t_{y_k}1_{v(g)\m})&=\af_{v(g)}(t_{y_k}1_{v(g)\m})1_{y_k}\\[.2em]
		&=\af_{\ta_{y_k}\ta\m_{y_l}u(g)}(t_{y_k}1_{v(g)\m})1_{\ta_{y_k}\ta\m_{y_l}}\\
		&\stackrel{\eqref{eq-identities}}{=}\af_{\ta_{y_k}\ta\m_{y_l}}\big(\af_{u(g)}(t_{y_k}1_{u(g)\m})1_{y_l}\big)1_{y_k},
	\end{align*}
	which implies that 
	\begin{align}\label{eq-fim-dem}
		\af_{(\ta_{y_k}\ta\m_{y_l})\m}\big(\af_{v(g)}(t_{y_k}1_{v(g)\m})\big)=\af_{u(g)}(t_{y_k}1_{u(g)\m}).
	\end{align}
	Observe that \eqref{eq-fim-dem} is also true for $h$ instead of $g$. Then,  \eqref{eq-strong} implies that
	\begin{align}\label{onemore}
		\af_{u(g)}(t_{y_k}1_{u(g)\m})f\neq 	\af_{u(h)}(t_{y_k}1_{u(h)\m})f.
	\end{align}
	Using that $t_{y_k}=\af_{\ta\m_{k,y}}(t_y)$ it is straightforward to check that \eqref{onemore} is equivalent to
	\begin{align}\label{onemoremore}
		\af_{\ta\m_{l,z}g}(t_y1_{(\ta\m_{l,z}g)\m})f\neq \af_{\ta\m_{l,z}h}(t_y1_{(\ta\m_{l,z}h)\m})f.	
	\end{align}
	Finally, since $\af_{\ta\m_{l,z}g}(t_y1_{(\ta\m_{l,z}g)\m})=\af_{\ta\m_{l,z}}\big(\af_g(t_y1_{g\m})\big)$ (and the same with $h$ instead of $g$), it follows from \eqref{onemoremore} that $\af_g(t_y1_{g\m})e\neq \af_h(t_y1_{h\m})e$. Consequently $T_y$ is $\af_{\G(y,z)}$-strong. This imply that $T$ is $\af$-strong, because $y,z \in \G_0$ were taken arbitrarily.
\end{proof}

\smallbreak

In the next two propositions we establish the Galois correspondence.

\begin{prop1} \label{correspondence-1} 
	Let $S$ be an $\alpha_{\G}$-partial Galois extension of $R=S^{\af_{\G}}$, $\cH\in \wSubg{\G}$ and $T=S^{\af_{\cH}}$. Then 
	\begin{enumerate}[\rm (i)]
		\item $T$ is $R$-separable and $\af$-strong,\smallbreak
		\item $\G_T=\cH$.
	\end{enumerate}
\end{prop1}
\begin{proof}
	
	\noindent By Lemma \ref{obs-Galois-compon}, $S_{y_j}$ is an $\af_{\G(y_j)}$-partial Galois extension of $R_{y_j}=S_{y_j}^{\af_{\G(y_j)}}$, for all $1\leq j\leq r$. Thus, \cite[Theorem 5.2]{DFP} implies that $T_{y_j}=S_{y_j}^{\af_{\cH_j(y_j)}}$ is  $R_{y_j}$-separable and $\af_{\G(y_j,y_j)}$-strong and $\G(y_j)_{T_{y_j}}=\cH_j(y_j)$, for all $1\leq j\leq r$. \smallbreak 
	
	Notice that (i) follows directly from  Lemma \ref{lem-separavel} and Proposition \ref{prop-strong}. The item (ii) follows from Proposition \ref{prop-importante} (ii).
\end{proof}

\begin{prop1} \label{correspondence-2} 
	Let $S$ be an $\alpha_G$-partial Galois extension of $R:=S^{\af_{\G}}$, $T$ an $R$-separable and $\af$-strong subring of $S$ such that $\G_T=\cH$, where ${\cH}\in \wSubg{\G}$. Then  $S^{\af_{\cH}}=T$.
\end{prop1}
\begin{proof}
	Let $y\in \G_0$, $T_y=T1_y$	and $R_y=R1_y$. Since $T$ is $\af$-strong, we have that $T_y$ is $\af_{\G(y)}$-strong. Also, $T_y$ is $R_y$-separable. In fact, consider $e=\sum_{i=1}^{m}a_i\otimes b_i\in T\otimes_R T$ an idempotent of separability of $T$ over $R$. It is immediate to verify that $e_y=\sum_{1\leq i\leq	m} (a_i)1_y\otimes (b_i)1_y \in T_y\otimes_{R_y} T_y$ is an idempotent of separability of $T_y$ over $R_y$. We claim that $G(y)_{T_y}=\cH_k(y)$, where $k\in \{1, 2, \ldots, r\}$ is taken such that $y\in Y_k$. It is clear that $\cH_k(y)\subset G(y)_{T_y}$. For the reverse inclusion, consider $g\in G(y)_{T_y}$ and $t\in T$. Notice that $\af_g(t1_{g\m})=\af_g(t(1_y1_{g\m}))=\af_g((t1_y)1_{g\m})=(t1_y)1_g=t(1_y1_g)=t1_g$.  
	Thus $T_y$ is $R_y$-separable and $\G(y,y)$-strong and $G(y)_{T_y}=\cH_k(y)$. Also, by Lemma \ref{obs-Galois-compon} we have that  $S_y$ is an $\af_{\G(y)}$-partial Galois extension of $S_y^{\G(y)}$. Hence, Theorem 5.4 of \cite{DFP} implies that $S_y^{\af_{\cH_k(y)}}=T_y$. On other hand, from Corollary \ref{prop-invariants-first} we obtain that 
	\[S^{\af_{\cH}}=\mathsmaller{ \bigoplus_{j=1}^r}\, S_j^{\af_{\cH_{j}}}, \qquad\text{where }\, S_j=\oplus_{z\in Y_j} S_z. \]
	By  Remark \ref{not-depend-choose-object}, $S_z^{\af_{\cH_k(z)}}=\af_{\ta_{k,z}}\big(\af_{\tau\m_{k,y}}\big(S_{y}^{\af_{\cH_k(y)}}\big)\big)$, for all $z\in Y_k$. Hence, we have that  
	$T_z=\af_{\ta_{k,z}}\big(\af_{\tau\m_{k,y}}\big(T_y\big)\big)$. From Corollary \ref{cor-invariants} follows that, $S_k^{\af_{\cH_{k}}}\subset \oplus_{z\in Y_k}\, T_z$, for all $1\leq k\leq r$. Consequently, $S^{\af_{\cH}}\subset T=\oplus_{y\in \G_0} T_y$. The converse inclusion is trivial.
\end{proof}

Let $R=S^{\af_{\G}}$ and denote by $\mathfrak{B}(S)$ the set of all subrings $T$ of $S$ which are $R$-separable, $\af$-strong and such that $\G_T=\cH$, for some $\cH\in \wSubg{\G}$.

\begin{teo1}\label{Galois-correspondence} {\rm \textbf{(Galois Correspondence)}} Let $S$ be an $\alpha_{\G}$-partial Galois extension of $S^{\af_{\G}}$. There exists a bijective correspondence  between $\wSubg{\G}$ and  $\mathfrak{B}(S)$ given by ${\cH}\mapsto S^{\af_{\cH}}$ whose inverse is given by $T\mapsto \G_T$.
\end{teo1}
\begin{proof}
The result follows directly from Propositions \ref{correspondence-1} and \ref{correspondence-2}.	
\end{proof}

\begin{remark}\label{gen-case}  {\rm \textbf{(General Case)}} {\rm The above theorem is also true when $\G$ is not connected. In fact,  suppose that $\G$ is not connected and $\af$ is a group-type partial action of $\G$ on $S$ in the sense of Definition \ref{def-groupoid-tipogrupo}. Assume also that $S$ is an $\alpha_{\G}$-partial Galois extension of $R=S^{\af_{\G}}$. In order to obtain Theorem \ref{Galois-correspondence} in this general case, we take the decomposition $\G=\dot\cup_{Z_j\in \G_0/\!\sim}\G_{Z_j}$ of $\G$ in its connected components, where $\G_0=Z_1\,\dot\cup\ldots \dot\cup\, Z_r$, and we fix: $\G_j:=\G_{Z_j}$, $S_j:=\oplus_{y\in Z_j}S_y$, $\af_j=(S_g,\af_{g})_{g\in \G_{j}}$ the partial action of $\G_{j}$ on $S_j$ and $R_j=S_j^{\af_j}$, for each $1\leq j\leq r$. By Proposition \ref{prop-reduction}, $S_j$ is an $\af_j$-partial Galois extension of $R_j$. Then, for a wide subgroupoid $\cH$ of $\G$ such that $\af_{\cH}$ is group-type, we have that $\cH_j=\cH\cap \G_j\in \wSubg{\G_j}$,  for all $1\leq j\leq r$. From Theorem \ref{Galois-correspondence} and Proposition \ref{prop-reduction}  follow that $T_j:=S_j^{\af_{\cH_j}}\in \mathfrak{B}(S_j)$. It is clear that $s=\sum_{j=1}^r s_j\in S^{\af_{\cH}}$ if and only if $s_j\in T_j$, for all $1\leq j\leq r$. Thus, $T=S^{\af_{\cH}}=T_1\oplus\ldots\oplus T_r$. Since $T_j$ is $R_j$-separable and $\af_j$-strong for all $j$, we obtain that $T$ is $R$-separable and $\af$-strong.  It is clear that $\G_T={\textstyle  \bigcup_{j=1}^r} (\G_j)_{T_j}$. Using that $ (\G_j)_{T_j}=\cH_j$, we conclude that $\G_T=\cH$. Conversely, consider a subring  $T$ of $S$ that is $R$-separable and $\af$-strong and $\G_T=\cH$, for some wide subgroupoid $\cH$ of $\G$ such that $\af_{\cH}$ is group-type. If  $T_j:=T1_{S_j}$, for all $1\leq j\leq r$, then $\G_{T_j}=\cH\cap \G_j\in \wSubg{\G_j}$. Hence	$T_j\in \mathfrak{B}(S_j)$. Thus, by Theorem \ref{Galois-correspondence},  $S^{\af_{\cH_j}}=T_j$, where $\cH_j=\G_{T_j}$. Therefore, $S^{\af_{\cH}}=\oplus_{j=1}^r S_j^{\af_{\cH_j}}=\oplus_{j=1}^r T_j=T$. } 
\end{remark}

Now we shall present several examples to illustrate the Galois correspondence given in the previous theorem.

\begin{exe}\label{exe1-trivial}{\rm
		Let $\G=\{x, y, g,g^{-1}\}$ be the groupoid with $\G_0=\{x,y\}$ and $g\in \G(x,y)$. Consider $S=ke_1\oplus ke_2$, where $k$ is a commutative ring, $e_ie_j=\delta_{i,j}e_i$ and $e_1+e_2=1$. Define the following partial action of $\G$ on $S$:
		\begin{align*}
			&S_{x}=ke_1=S_{g\m},& &S_{y}=ke_2=S_{g},& &\af_x=\id_{S_x},& &\af_y=\id_{S_y},& &\af_g(e_1)=e_2.&
		\end{align*}
		Thus, $R:=S^{\af_{\G}}=k(e_1+e_2)$ and $R\subset S$ is an $\alpha$-partial Galois extension. In fact, it is immediate to check that $a_i= b_i = e_i$, $i \in \{1, 2\}$, determine a partial Galois coordinate system of $S$ over $R$. Also, $\af$ is a group-type partial action with transversal $\tau=\{\tau_x=x,\tau_y=g\}$ for $x$. In this case, the Galois correspondence is given by 
		\begin{align*}
			&{{\cH}}=\G\longleftrightarrow {T}=k(e_1+e_2),& 	&{{\cH}}=\{x,y\}\longleftrightarrow {T}=S.&
	\end{align*}}
\end{exe}

\begin{exe}\label{exe2-global}{\rm
		Let $\G$ be the groupoid considered in Example \ref{ex:invariant} and consider $S:=\oplus_{i=1}^4k e_i$, where $k$ is a commutative ring,  $e_ie_j=\delta_{i,j}e_i$ and $\sum_{i=1}^{4}e_i=1$. Define the following (global) action $\beta$ of $\G$ on $S$:
		\begin{align*}
			&\beta_x=\id_{S_x},\,\,\,S_{x}=k e_1\oplus  k e_2,&  &\beta_y=\id_{S_y},\,\,\,S_{y}=k e_3\oplus k e_4, & \\
			&\beta_g: k e_1\oplus k e_2\to k e_1\oplus k e_2,& &\beta_g(e_1)=e_2,\,\,\,\beta_g(e_2)=e_1, & \\
			&\beta_h: k e_3\oplus k e_4\to k e_3\oplus k e_4,& &\beta_h(e_3)=e_4,\,\,\,\beta_g(e_4)=e_3,&\\
			&\beta_l: k e_1\oplus k e_2\to k e_3\oplus k e_4,& &\beta_l(e_i)=e_{i+2},\,\,\,i=1,2.&  
		\end{align*}
		Since $m=lg$, we put $\beta_m=\beta_l\beta_g$. Notice that $S^{\af_{\G}}=k$ and $k\subset S$ is an $\alpha$-partial Galois extension. In fact, it is easy to check that $a_i = b_i = e_i$, $1 \leq i \leq 4$, is a partial Galois coordinate system of $S$ over $k$. 
		Also, $\beta$ is a group-type partial action because $\beta$ is global. In this case, the Galois correspondence is given by 
		\begin{align*}
			{{\cH}}=\G&\longleftrightarrow {T}=k&\\[.2em]
			{{\cH}}=\{x,y,l,l\m\}&\longleftrightarrow {T}=k(e_1+e_3)\oplus k(e_2+e_4)&\\[.2em]
			{{\cH}}=\{x,y,m,m\m\}&\longleftrightarrow {T}=k(e_1+e_4)\oplus k(e_2+e_3)&\\[.2em]	
			{{\cH}}=\{x,y,g,h\}&\longleftrightarrow {T}=k(e_1+e_2)\oplus k(e_3+e_4)&\\[.2em]	
			{{\cH}}=\{x,y,g\}&\longleftrightarrow {T}=k(e_1+e_2)\oplus ke_3\oplus k e_4&\\[.2em]	
			{{\cH}}=\{x,y,h\}&\longleftrightarrow {T}=ke_1\oplus k e_2\oplus k(e_3+e_4) &\\[.2em]
			{{\cH}}=\{x,y\}&\longleftrightarrow {T}=S&
	\end{align*}}
\end{exe}

\begin{exe}\label{ex:groupoid-12}{\rm
		Let $\G=\{g,g^2,h,h^2,l,l\m,m,m\m,n,n\m\}$ be the groupoid with set of objects $\G_0=\{x,y\}$. The composition in $\G$ is given by
		\begin{align*}
			g&\in \G(x),\quad h\in \G(y),& &\,\,\,g^3=x,& &\,h^3=y,& \\
			&l,m,n\in\G(x,y),& lg&=m=hl,& mg&=n=hm.&
		\end{align*} 
		The groupoid $\G$ is illustrates in the following diagram 
		\[\xymatrix{& x\ar[d]^{n}\ar[r]^{g}  &x\ar[d]^{m} \ar[r]^{g} &x\ar[d]^{l}\\
			& y & y\ar[l]^{h}&\ar[l]^{h} y}\\ \]
		Let  $S:=\oplus_{i=1}^6k e_i$, where $k$ is a commutative ring, $e_ie_j=\delta_{i,j}e_i$ and $\sum_{i=1}^{6}e_i=1$. Define the following partial action $\alpha$ of $\G$ on $S$:
		\begin{align*}
			&\af_x=\id_{S_x},\,\,\,S_{x}=k e_1\oplus&  &\!\!\!\!k e_2\oplus  k e_3,&  &\af_y=\id_{S_y},\,\,\,S_{y}=k e_4\oplus k e_5\oplus  k e_6, & \\
			&\af_l: k e_1\oplus k e_2\oplus k e_3\to& &\!\!\!\!k e_4\oplus k e_5\oplus k e_6,& &\af_l(e_i)=e_{i+3},\,\,\,i=1,2,3.& \\
			&\af_g: k e_1\to k e_2,& &\!\!\!\!\af_g(e_1)=e_2,& & \af_{g^2}=\af\m_g,& \\
			&\af_h: k e_4\to k e_5,& &\!\!\!\!\af_h(e_4)=e_5,& &\af_{h^2}=\af\m_h,&\\
			&\af_m: k e_1\to k e_5,& &\!\!\!\!\af_m(e_1)=e_5,& &\af_{m\m}=\af\m_m&\\
			&\af_n: k e_2\to k e_4,& &\!\!\!\!\af_n(e_2)=e_4,& &\af_{n\m}=\af\m_n.&
		\end{align*}
		Notice that $R=S^{\af_{\G}}=k(e_1+e_2+e_4+e_5)\oplus k(e_3+e_6)$ and $R\subset S$ is an $\alpha$-partial Galois extension. In fact, it is immediate to verify that $a_i = b_i = e_i$, $1 \leq i \leq 6$, is a partial Galois coordinate system of $S$ over $R$. 
		Also, $\af$ is a group-type partial action because $\tau(x)=\{\tau_x=x,\tau_y=l\}$ is a transversal in $\G$ for $x$ that satisfies \eqref{cond1}. The Galois correspondence is given by: 
		\begin{align*}
			\G&\longleftrightarrow k(e_1+e_2+e_4+e_5)\oplus k(e_3+e_6)&\\[.2em]
			\{x,y,l,l\m\}&\longleftrightarrow k(e_1+e_4)\oplus k(e_2+e_5)\oplus k(e_3+e_6)&\\[.2em]
			\{x,y,g,g^2,h,h^2\}&\longleftrightarrow k(e_1+e_2)\oplus k e_3\oplus k(e_4+e_5)\oplus k e_6&\\[.2em]
			\{x,y,g,g^2\}&\longleftrightarrow k(e_1+e_2)\oplus k e_3\oplus k e_4\oplus k e_5\oplus k e_6&\\[.2em]
			\{x,y,h,h^2\}&\longleftrightarrow ke_1\oplus ke_2\oplus k e_3\oplus k(e_4+e_5)\oplus k e_6& \\[.2em]
			\{x,y\}&\longleftrightarrow S&
		\end{align*}
		Observe that $\cH_1=\{x,y,m,m\m\}$ and $\cH_{2}=\{x,y,n,n\m\}$ are wide subgroupoids of $\G$. However $\af_{\cH_{1}}$ and $\af_{{\cH}_2}$ are not group-type partial actions. Consequently $\cH_{1}, \cH_{2}\notin \wSubg{\G}$.}
\end{exe}

\begin{exe}\label{ex_inv_semigroup} {\rm (Based on \cite[Example 5.15]{LT})
		Let $I(X)$ be the groupoid of the partial bijections of $X = \{1,2,3\}$. Let $\G$ be the subgroupoid of $I(X)$ given by the bijections between subsets of two elements of $X$. The composition of two maps $f$ and $g$ is defined in $\G$ when the domain of $f$ is equal to the image of $g$. We fix $X_{i,j}:=\{i,j\}$ for all $1\leq i< j\leq 3$. The identity map on $X_{i,j}$ is denoted by $\id_{ij}$ and we define the following bijections:
		\begin{align*}
			&f_{ij} : X_{i,j} \rightarrow X_{i,j},\quad\,\,\, f_{ij}(i)=j,\quad\,\,\,\, f_{ij}(j)=i,\quad\,\,\,\text{for all } 1\leq i< j\leq 3,	\\
			&d_{ik}^j : X_{i,j} \rightarrow X_{j,k},\quad d_{ik}^j(i)=k,\quad d_{ik}^j(j)=j,\qquad\text{for all } i\neq j,\,i\neq k, j\neq k,\\
			&p_{ij}^{kl} : X_{i,j} \rightarrow X_{j,k},\quad p_{ij}^{kl}(i)=k,\quad p_{ij}^{kl}(j)=l,\qquad\, \text{where }  X_{i,j}\cap X_{k,l}=\{i\}\text{ or } \{j\}.
		\end{align*}
		Denote $x=\id_{12}$, $y=\id_{13}$ and $z=\id_{23}$. The set of objects of $\G$ is $\G_0=\{x,y,z\}$ and we have the following morphisms:
		\begin{align*}
			&f_{12}\in \G(x),& &f_{13}\in \G(y),& &f_{23}\in \G(z),& &d_{13}^2\in \G(x,z),& &d_{23}^1\in \G(x,y),&\\ 
			&d_{12}^3\in \G(y,z),&  &d_{31}^2\in \G(z,x),& &d_{32}^1\in \G(y,x),& &d_{23}^1\in \G(z,y),& &p_{12}^{13}\in \G(x,y),&\\
			&p_{12}^{23}\in \G(x,z),& &p_{13}^{23}\in \G(y,z),& &p_{13}^{12}\in \G(y,x),& &p_{23}^{12}\in \G(z,x),& &p_{23}^{13}\in \G(z,y).&
		\end{align*}
		Notice that $\G$ is connected and $\G \simeq \G_0^2 \times C_2$, where $C_2$ denotes the cyclic group of order $2$. We will construct a partial action $\af$ of $\G$ on the ring $S = \oplus_{i = 1}^6 ke_i$, where $k$ is a commutative ring, $e_ie_j = \delta_{i,j}e_i$ and $\sum_{i = 1}^6 e_i = 1_S$. Define $\af_x=\id_{S_x}$, $\af_y=\id_{S_y}$ and $\af_z=\id_{S_z}$, where
		\[S_{x} = ke_1 \oplus ke_2,\qquad S_{y} = ke_3 \oplus ke_4,\qquad S_{z} = ke_5 \oplus ke_6.\] 
		The partial isomorphism $\af_{f_{ij}}$ are defined by:
		\begin{align*}
			&\alpha_{f_{12}} : ke_1\oplus ke_2 \rightarrow ke_1\oplus ke_2,\qquad e_1\mapsto e_{2},\quad e_2\mapsto e_{1},&  \\[.2em]
			&\alpha_{f_{13}} : ke_3\oplus ke_4 \rightarrow ke_3\oplus ke_4,\qquad e_3\mapsto e_{4},\quad e_4\mapsto e_{3},&  \\[.2em]
			&\alpha_{f_{23}} : ke_5\oplus ke_6 \rightarrow ke_5\oplus ke_6,\qquad e_5\mapsto e_{6},\quad e_6\mapsto e_{5}.& 
		\end{align*}
		The partial isomorphism $\af_{d_{ijk}}$ are defined by:
		\begin{align*}
			&\alpha_{d_{13}^2} : ke_1 \oplus ke_2 \rightarrow ke_5  \oplus ke_6, \quad e_1\mapsto e_{6},\,\,e_2\mapsto e_5,&  &\af_{d_{31}^2}=\big(\alpha_{d_{13}^2}\big)\m,&\\[.2em]
			&\alpha_{d_{23}^1} : ke_1 \oplus ke_2 \rightarrow ke_3  \oplus ke_4, \quad e_1\mapsto e_{3},\,\,e_2\mapsto e_4,&  &\af_{d_{32}^1}=\big(\alpha_{d_{23}^1}\big)\m,&\\[.2em]
			&\alpha_{d_{12}^3} : ke_3 \rightarrow ke_5,\qquad\qquad\quad\,\,\, \quad e_3\mapsto e_{5},&  &\af_{d_{21}^3}=\big(\alpha_{d_{12}^3}\big)\m.&
		\end{align*}
		The partial isomorphism $\af_{p_{ij}^{kl}}$ are defined by:
		\begin{align*}
			&\alpha_{p_{12}^{13}} : ke_2 \rightarrow ke_3,\quad\quad e_2\mapsto e_{3},&  &\af_{p_{13}^{12}}=\big(\alpha_{p_{12}^{13}}\big)\m,&\\[.2em]
			&\alpha_{p_{12}^{23}} : ke_1 \rightarrow ke_5,\quad\quad e_1\mapsto e_{5},&  &\af_{p_{23}^{12}}=\big(\alpha_{p_{12}^{23}}\big)\m,&\\[.2em]
			&\alpha_{p_{13}^{23}} : ke_3 \rightarrow ke_6,\quad\quad e_3\mapsto e_{6},&  &\af_{p_{23}^{13}}=\big(\alpha_{p_{13}^{23}}\big)\m,&
		\end{align*}
		Observe that $S^{\af_G} =k$ and $k \subset S$ is an $\alpha_{\G}$-partial Galois extension. In fact, it is easy to check that $a_i = b_i = e_i$, $1 \leq i \leq 6$, is a partial Galois coordinate system of $S$ over $k$. Also, $\tau(x) = \{\tau_{x}=x, \tau_{y} = d_{13}^2, \tau_{z} = d_{23}^1\}$ is a transversal in $\G$ for $x$ that satisfies \eqref{cond1}, that is, $\af$ is a group-type partial action. We denote $e_{i_1,\ldots,i_n}=e_{i_1}+\ldots+e_{i_n}$. The Galois correspondence is given by: 
		\begin{align*}
			\cH_0=\G_0&\longleftrightarrow T_0=S&\\[.2em]
			\cH_1=\G_0\cup \{ f_{12} \}&\longleftrightarrow T_1=ke_{1,2} \oplus ke_3 \oplus ke_4  \oplus ke_5  \oplus ke_6&\\[.2em]
			\cH_2=\G_0\cup \{ f_{13} \}&\longleftrightarrow T_2=ke_{3,4}  \oplus ke_1 \oplus ke_2  \oplus ke_5  \oplus ke_6&\\[.2em]
			\cH_3=\G_0\cup \{ f_{23} \}&\longleftrightarrow T_3=ke_{5,6} \oplus ke_1 \oplus ke_2  \oplus ke_3  \oplus ke_4&\\[.2em]
			\cH_4=\G_0\cup \{d_{23}^1, d_{32}^1 \}&\longleftrightarrow T_4=ke_{1,3}  \oplus ke_{2,4}\oplus ke_5 \oplus ke_6&\\[.2em]
			\cH_5=\G_0\cup \{d_{13}^2, d_{31}^2 \}&\longleftrightarrow T_5=ke_{1,6} \oplus ke_{2,5}\oplus ke_3 \oplus ke_4&\\[.2em]
			\cH_{6}=\cH_1\cup \cH_2&\longleftrightarrow T_{6}=ke_{1,2} \oplus ke_{3,4}\oplus ke_5 \oplus ke_6&\\[.2em]
			\cH_{7}=\cH_1\cup \cH_3&\longleftrightarrow T_{7}=ke_{1,2}\oplus ke_{5,6}\oplus ke_3 \oplus ke_4&\\[.2em]
			\cH_{8}=\cH_2\cup \cH_3&\longleftrightarrow T_{8}=ke_{3,4} \oplus ke_{5,6}\oplus ke_1 \oplus ke_2&\\[.2em]
			\cH_{9}=\cH_2\cup \cH_5&\longleftrightarrow T_{9}=ke_{1,6}\oplus ke_{3,4}\oplus ke_2 \oplus ke_5&\\[.2em]
			\cH_{10}=\cH_3\cup \cH_4&\longleftrightarrow T_{10}=ke_{1,3}\oplus ke_{5,6}\oplus ke_2 \oplus ke_4&\\[.2em]
			\cH_{11}=\cH_4\cup \cH_5\cup \{ p_{13}^{23}, p_{23}^{13}  \} &\longleftrightarrow T_{11}=ke_{1,3,6} \oplus ke_{2,4,5}&\\[.2em]
			\cH_{12}= \cH_1\cup \cH_2 \cup \cH_4\cup \{ p_{12}^{13}, p_{13}^{12}  \}  &\longleftrightarrow T_{12}=ke_{1,2,3,4} \oplus ke_5 \oplus ke_6&\\[.2em]
			\cH_{13}=  \cH_1\cup \cH_3 \cup  \cH_5\cup  \{ p_{12}^{23}, p_{23}^{12}  \}  &\longleftrightarrow T_{13}=ke_{1,2,5,6} \oplus ke_3 \oplus ke_4&\\[.2em]
			\cH_{14}= \cH_2\cup \cH_3 \cup \{d_{12}^3, d_{21}^3, p_{13}^{23}, p_{23}^{13}  \}   &\longleftrightarrow T_{14}=ke_{3,4,5,6} \oplus ke_1 \oplus ke_2&\\[.2em]
			\cH_{15}= \cH_1 \cup \cH_2\cup \cH_3 & \longleftrightarrow T_{15}=ke_{1,2} \oplus ke_{3,4} \oplus ke_{5,6}&\\[.2em]
			\cH_{16}=  \cH_4 \cup \cH_{15} \cup  \{ p_{12}^{13}, p_{13}^{12}  \}  &\longleftrightarrow T_{16}=ke_{1,2,3,4} \oplus ke_{5,6}&\\[.2em]
			\cH_{17}=  \cH_5 \cup \cH_{15} \cup  \{ p_{12}^{23}, p_{23}^{12}  \}  &\longleftrightarrow T_{17}=ke_{1,2,5,6} \oplus ke_{3,4}&\\[.2em]
			\cH_{18}=  \cH_{15} \cup \{d_{12}^3, d_{21}^3, p_{13}^{23}, p_{23}^{13}  \}  &\longleftrightarrow T_{18}=ke_{3,4,5,6} \oplus ke_{1,2}&\\[.2em]
			\cH_{19}=\G&\longleftrightarrow T_{19}=k.&
	\end{align*}}
\end{exe}

\subsection{Application to the global case} In this subsection we will assume that $\alpha$ is a global action of a non-necessarily connected groupoid $\G$ on $S$. By Lemma 1.1 of \cite{BP}, $S_g=S_{t(g)}$, for all $g\in \G$. In particular, $\alpha$ and $\alpha_{\cH}$ are group-type, for all subgroupoid $\cH$ of $\G$. Thus $\cH\in\wSubg{\G}$ if and only if $\cH$ is wide. Also, by Proposition \ref{prop-strong}, $S^{\af_{\cH}}$ is $\af$-strong in the sense that was defined in p. 108 of \cite{PT}. On the other hand, consider a subring $T$ of $S$ which is $S^{\alpha_{\G}}$-separable and $\af$-strong. By Proposition \ref{prop-importante} (i) we have that $\G_T$ is a wide subgroupoid of $\G$. Using Remark \ref{gen-case} we have the following result. 

\begin{prop1} \label{global-case-1} Suppose that $S$ is an $\alpha_{\G}$-Galois extension of $R=S^{\alpha_{\G}}$. The associations $\cH \mapsto S^{\af_{\cH}}$ and $T\mapsto \G_{T}$ determine a bijective correspondence between the wide subgroupoids of $\G$ and the subrings of $S$ that are $\alpha$-strong and $R$-separable.  
\end{prop1}

The correspondence given in the above proposition was proved in Theorem 4.6 (i) of \cite{PT}. However, using our results, we can determine the elements $S^{\af_{\cH}}$ and $\G_{T}$ of this correspondence from the classical Galois theory of groups given in \cite{CHR}. \vspace{.05cm}

Assume that $S$ is an $\alpha_{\G}$-Galois extension of $R=S^{\alpha_{\G}}$ and  $\G=\dot\cup_{Z_j\in \G_0/\!\sim}\G_{Z_j}$ is the decomposition of $\G$ in connected components, where $\G_0=Z_1\,\dot\cup\ldots \dot\cup\, Z_r$. Consider $\G_j:=\G_{Z_j}$, $S_j:=\oplus_{y\in Z_j}S_y$, $\af_j=(S_g,\af_{g})_{g\in \G_{j}}$ the partial action of $\G_{j}$ on $S_j$ and $R_j=S_j^{\af_j}$, for each $1\leq j\leq r$. Let $\cH$ be a wide subgroupoid of $\G$. In this case, $\af_{\cH}$ is group-type. By Remark \ref{gen-case}, $T=S^{\af_{\cH}}=T_1\oplus\ldots\oplus T_r$ with $T_j=S_j^{\af_{\cH_j}}$ and $\cH_j=\cH\cap \G_j$. If $\left(\cH_j\right)_0=Y_{j,1}\,\dot\cup\,\ldots\,\dot\cup\, Y_{j,n_j}$ then the groupoid $\cH_j$ admits a decomposition $\cH_{j}=\cH_{j,1}\dot\cup\,\ldots\,\dot\cup\cH_{j,n_j}$ in connected components, for all $1\leq j\leq r$. Let $y_{j,i}\in Y_{j,i}$ for all $1\leq i\leq n_j$. By Corollary \ref{prop-invariants-first}, \[T_j=S_j^{\af_{\cH_j}}\simeq {\textstyle \bigoplus_{i=1}^{n_j}} S_{y_{j,i}}^{\alpha_{\cH_{j,i}(y_{j,i})}}.\] Observe that the group $\cH_{j,i}(y_{j,i})$ acts globally on $S_{y_{j,i}}$. Thus, $S^{\af_{\cH}}$ is completely determined by invariant elements of global group actions. In fact, we have the following. \vspace{0.05cm}

\begin{prop1} \label{global-case-2} Assume the above notation. Then 
	\[S^{\af_{\cH}}\simeq  \mathsmaller{ \bigoplus_{j=1}^{r} \bigoplus_{i=1}^{n_j}} S_{y_{j,i}}^{\alpha_{\cH_{j,i}(y_{j,i})}}. \]
\end{prop1}
\vspace{0.05cm}

Let $\tau_{j,i}=\{\ta_{j,i,y}:y_{j,i}\to y\}_{y\in Y_{j,i}}$ be a transversal for $y_{j,i}$ in $\cH_{j,i}$ for all $1\leq i\leq n_j$ and $1\leq j\leq r$. Consider $T$ a subring of $S$ which is $S^{\af_{\G}}$-separable and $\af$-strong and consider $T_j=T1_{S_j}$. From Lemma \ref{lem-invari-subgroupoid} it follows that \[(\G_j)_{T_j}= {\textstyle    \bigcup_{i=1}^{n_j} }\, {\textstyle \bigcup_{y\in Y_{j,i}}}\, \ta_{j,i,y}\G_j(y_{j,i})_{(T_j)_{y_{j,i}}}\ta\m_{j,i,y}.\] Observe that $\G_j(y_{j,i})_{(T_j)_{y_{j,i}}}$ is the subgroup of $\G_j$ whose elements fix $(T_j)_{y_{j,i}}$. Since $\G_T={\textstyle    \bigcup_{j=1}^{r} }\,(\G_j)_{T_j}$, we obtain the following.
\vspace{0.05cm}

\begin{prop1} \label{global-case-3} Assume the above notation. Then 
	\[\G_T= {\textstyle    \bigcup_{j=1}^{r} } \,{\textstyle    \bigcup_{i=1}^{n_j} }\, {\textstyle \bigcup_{y\in Y_{j,i}}}\, \ta_{j,i,y}\G_j(y_{j,i})_{(T_j)_{y_{j,i}}}\ta\m_{j,i,y}.\]
\end{prop1}

We finish the paper noting that the problem of determining a Galois correspondence for partial groupoid actions (not necessarily of group-type) remains an open question.

\end{document}